\RequirePackage{fix-cm}
\documentclass[preprint]{imsart}       
\usepackage{graphicx}
\usepackage[margin=1.5in]{geometry}
\usepackage{dsfont}
\usepackage{mathtools}
\usepackage{amssymb}
 \usepackage{mathptmx}      
 \usepackage{amsthm}
  \usepackage{url}
\usepackage{hyperref}
\newtheorem{theorem}{Theorem}[section]
\newtheorem{definition}[theorem]{Definition}
\newtheorem{proposition}[theorem]{Proposition}

\newtheorem{lemma}[theorem]{Lemma}
\newtheorem{corollary}[theorem]{Corollary}
\newtheorem{remark}[theorem]{Remark}
\begin{document}

\title{Fluid limits for Queue-based CSMA with polynomial rates, homogenization and reflection} 

\runtitle{\large Fluid limits of QB-CSMA}
\begin{aug}
	\author[A]{\fnms{Eyal} \snm{Castiel}\ead[label=e1]{eyal.ca@campus.technion.ac.il}},
	\address[A]{ISAE Supaero, 10 Avenue Edouard Belin, 31400 Toulouse, France
	}
	
	\address[B]{Universite de Toulouse Paul Sabatier, 118 Route de Narbonne, 31400 Toulouse, France
		\printead{e1}}
\end{aug}




\maketitle

\begin{abstract}	\large
	We study in this paper a variation of the acclaimed CSMA random access protocol. We will focus on the case where back-off rates at each node is polynomial in the size of the queue. Under a condition relating the exponent in the polynomial rates and the geometry of the interference graph, we prove convergence of the scaled process to a deterministic fluid limit up to the time a queue reaches $0$ on the fluid scale . We outline the difficulties arising at that time and solve them in the case of a complete interference graph. This paper relies on a new method to obtain a fully coupled stochastic averaging principle and can hopefully lead to more result in heavy load situations.

\end{abstract}
\renewcommand{\d}{{\rm{d}}}

\newcommand{\ind}{\mathds{1}}
\newcommand{\reels}{\mathds{R}}
\newcommand{\entier}{\mathds{N}}
\newcommand{\cesp}[2]{\mathds{E}_{#1}\left[#2\right]}
\newcommand{\norm}[2]{\left\lVert #1\right\rVert_{#2}}
\newcommand{\cpro}[2]{\mathds{P}_{#1}\left(#2\right)}
\renewcommand{\P}{\mathds{P}}
\newcommand{\E}{\mathds{E}}

	
	{\large
	\setcounter{tocdepth}{1}
	\tableofcontents
	
	\section{\large{Introduction}}
	In this paper, we will prove fluid limit results for the QB-CSMA introduced in \cite{adia}, with polynomial rates using the homogenization result from \cite{thesis}, a special case has already been investigated in Heavy traffic using the same method in \cite{us}. In this paper, we are concerned with a scheduling problem where users are situated on the nodes of a simple non-oriented graph. Each node of the graph has a stream of Poisson arrivals of jobs that require a service time exponential with parameter one. A node (also called server) is said to be active if it is processing jobs. Edges in the graph signify the impossibility of simultaneous activation on both nodes sharing the edge. One of the main application of this model is the study of wireless networks. Once a queue start transmitting it prevents all of its neighbours from starting a transmission. One algorithm still widely used today to schedule conflicting queues on a graph are CSMA algorithms. Each node use carrier sensing abilities to listen to availability of the channel and has a \textit{fugacity} associated which dictates how often and how long it activates to transmit packets. The type of CSMA algorithm studied in this paper is \textit{queue-based}, meaning the activation and deactivation rate at each node depends on the current queue length at that node. The main difficulty of the results presented consists in establishing a homogenization result when queue lengths are away from zero and bridging the behavior of the process at and away from zero. We will use a previously established homogenization result and accomplish the second goal with a careful coupling and a tightness argument. 
	\subsection{\underline{\large{Motivation and context}}}
	The search for efficient scheduling algorithm with provable performance bounds had a major breakthrough with the introduction of the Max-Weights algorithm by Tassiulas and Ephremides in \cite{mw}. However the centralized nature of this algorithm and its reliance on solving a potentially NP-hard optimization problem at each iteration makes its implementation in large scale networks challenging. One of the issue of classical CSMA algorithm is that it is not possible to design fugacities that ensure queue lengths are stable for any subcritical arrival rates , even if the graph topology is known. For nearly twenty years, the search for distributed algorithm with maximum throughput capacity has been unsuccessful until Jiang and Walrand \cite{jiang10} and Rajagopolan et al \cite{adia}. Their idea was to approximate the decisions made by Max-Weight using only local information of the local queue and the transmitting state of neighbours. 
	
	The operations of distributed scheduling algorithms in these scenarios
	may be described as follows.
	At any time, a node either starts a random
	back-off period or proceeds with the transmission of the current packet,
	if any, with a probability that depends on the local queue length. The order of service of packets at each queue is \textit{First come first served} (FCFS).
	When inactive, a node simply runs down its back-off clock, but freezes
	it whenever any of its neighbours in the conflict graph are active,
	ensuring that a back-off period can only end when all its neighbours
	are inactive.
	At that point, a node either initiates a packet transmission
	or proceeds with the next random back-off period with a probability
	which is also a function of the local queue length.
	
	Now observe that transmission activity in wireless networks can be
	detected by `sensing' a shared channel, and that the above back-off
	mechanism precludes concurrent transmissions of mutually interfering nodes.
	Further note that the idleness and randomized deactivation may seem
	inefficient from a resource utilization perspective, but play
	an instrumental role in sharing the medium through `listening'
	in the absence of any centralized access control mechanism. The extreme case where a node deactivates only when its queue is empty is referred to as the Random Capture algorithm~\cite{feurca}. Although it may seem to minimize idleness, it was actually shown that it is not always throughput-optimal depending on the graph topology. The authors of ~\cite{Ghaderi2012} exhibit such a topology with a ''broken diamond`` network.
	
	Beside the throughput, little is known about the queueing dynamic of this algorithm. This paper aims at deepening our understanding of this promising algorithm. Two situations will be studied: a general and a full interference graphs. The main insight in this paper is that the homogenization technique used here are extremely efficient as long as queue lengths are bounded away from zero. With minimal information on the interference graph (for instance number of nodes or maximum size of an independent set), starting from a positive initial condition at each node, it is easy to obtain convergence of the scaled process to a deterministic function governed by an ODE up to the time the solution hits $0$ in one of its coordinates. In the case of a complete interference graph, we are able to go further and distinguish (sub/super)critical cases. Each case gives a different limiting behaviour. We give a proof for convergence to a fluid limit for any initial condition, any arrival rates and any finite time horizon. Possible reflections at zero are the main consideration of this paper and are handled with an ad-hoc coupling. It is very likely that under some ``worst case scenario'' conditions on the arrival rates, similar reasoning will lead to a result for heavily loaded systems in general graph topologies. Similarly to \cite{us}, we expect a non-trivial behavior at a time scale intermediate between $N$ and $N^2$. In this paper, the process of interest will be denoted $(Q^N,\sigma^N)$ and is given by 
	\[Q^N(t)\coloneqq\frac{Q(Nt)}{N}\text{ and } \sigma^N(t)\coloneqq \sigma(Nt).\]
	To be rigorous, we would have to add a superscript to the unscaled process: we actually consider a sequence of networks where the initial condition may depend on $N$ the scaling parameter. For the sake of clarity of exposition, we will omit this dependency in $N$ for the initial state. We will consider situations where $Q^N(0)=\frac{Q(0)}{N}\to q^0\geqslant 0$ as $N\to +\infty$.

	\subsubsection{\large{\underline{Queue-Based CSMA and polynomial rates}}}
The algorithm of Queue-Based CSMA has been introduced in \cite{adia} and tried to answer the question of finding a distributed algorithm approximating the Max-Weights algorithm. The main concern at the time was finding a distributed algorithm that has the largest stability region possible. This was proved in \cite{shashin}, \cite{ghad10} for increasingly ``aggressive'' activity functions. In this context, aggressive means rapidly increasing at infinity. It has been empirically tested and proved in steady state in some cases \cite{Bouman11} that a more rapidly growing activation function leads to better delay bounds up to a certain point. It is commonly conjectured that this should be true in any topology. In some graph topology, the activity function being too agressive can lead to a loss in the stability region. A rapidly growing activity function can lead to starvation of some queues and even instability \cite{Ghaderi2012}. One way to prevent starvation is to ensure that the scheduler moves ``fast'' between possible activity states. In those situations, we can study the system as if the scheduler interacts with queue lengths through a steady state average (see the next section for more details). When the activity function grows faster than $f(x)=x$ a queue will be almost empty by the time it deactivates so polynomial rates $f(x)=x^a$, $a<1$ are the fastest growing activity functions we can hope for while maintaining this property. The question of throughput optimaly for CSMA-type algorithm has already been widely studied, see for instance \cite{Yun12:0} for a survey of the subject. The ``time scale separation'' assumption is here justified rigorously through homogenization.
	\subsection{\large{\underline{Polling}}}
	
	When the interference graph is a complete graph, at most one queue can be active at any given time. This is also called a polling model. There are many ways to schedule and distribute the available resource betweeen the competing nodes: exhaustive service (the server only releases the channel when it is empty), gated service (the server releases the channel when all the jobs that were present when it started serving have been processed), fixed service (server processes a fixed number of jobs before releasing), threshold service (the server provides service until its backlog becomes smaller than a fixed threshold), the server may process a random number of requests, etc $\ldots$ When a queue deactivates, there may or may not be some idling time before the next queue activates. In the literature there is a distinction between polling models with and without switch-over times.  Polling systems are however often studied as centralized models with policies like SRPT (Shortest Remaining Processing Time), FCFS (First Come First Serve), and others... See \cite{polbox} for an overview of a Mean Value Analysis of different scheduling policy for polling systems.
	
	
	On a complete interference graph, CSMA and QB-CSMA are in fact a polling systems with switchover time and random activation duration. This equivalence has been exploited to use results for polling systems where the server only processes one job before moving to the next queue and a probabilistic routing policy. It has been used to analize CSMA algorithms where nodes deactivate at a fixed (non-queue-based) rate, see for instance~ \cite{Dorsman15}.
	\subsubsection{\large{\underline{Stochastic averaging principle}}}\label{sec:sap}
	The analysis of the system is significantly simplified by a stochastic averaging phenomenon. This type of behaviour can occur when two interacting processes operate on different time scale. This phenomenon is said to occur when the 'fast' process, in our case the scheduler, only interact with the 'slow' process, in our case the queue lengths, through an instantaneous equilibrium. The type of stochastic averaging studied here is of the most difficult kind because the dynamic of the fast process also depends on the value of the slow component. In a simple sense, for Markov processes, the idea is to replace the generator of the slow variable $L^\sigma_{\mathrm{s}}f(q)$ by a homogenized version: $\sum_{\sigma}\pi^q(\sigma)L^\sigma_{\mathrm{s}}f(q)$, and completely forget the fast variable.
	
	This type of approximation is in general difficult to obtain and explicit bounds are all the more valuable. Numerous methods have been developed to that effect with and without explicit bounds. We mention for instance the classical monograph by Freidlin and Wentzell \cite{freiwent}, the standardized method by Kurtz in \cite{Kurtz92} and the method yielding explicit bounds by Luczak and Norris \cite{Luc13}. In this paper, we will use the most recent method developed in \cite{us} and in more details in \cite{thesis}. This method, more tailored towards Markov processes than the one in \cite{Luc13}, uses solutions to Poisson equation and leads to explicit bounds expressed simply in terms quantities depending on the slow and fast dynamics. The focus of this paper is an application of this method with the added difficulty of dealing with the boundary. The intuition is that we need to check that the invariant measure of $L^{Q_t}_{\mathrm{f}}$ changes ``slowly'' compared to the mixing time of $L^{Q_t}_{\mathrm{f}}$, so that the schedule reaches equilibrium before the ``goal measure'' changes. \cite{us} and \cite{thesis} give an upper bound depending on the interference graph for the range of exponent such that QB-CSMA with activation function $f(x)=x^a$ is subject to a stochastic averaging principle when queue lengths are bounded away from zero. Unfortunately, the 'speed' at which the invariant measure of $L^{Q_t}_{\mathrm{f}}$ evolves depends on the minimum of queue lengths.  
	\subsubsection{\large{\underline{Fluid limits}}}
	When considering stability of queueing networks, in terms of tightness of its distribution, it is common to consider a rescaling in time and space of the process and let the scaling parameter go to $+\infty$. When the scaling is the same in time and space, the limiting points (if they exist) are called ``fluid limits''. When unique, the fluid limit can serve as a first order approximation of the process for large $N$ and stability of the original process can be deduced from the behaviour of the fluid limit. This method is more than 40 years old: it has been introduced in \cite{malyshev} where the authors use it to establish the transient/recurrent behaviour of random walks on $\mathds{Z}^2$ and $\mathds{Z}^3$.
	
	Fluid limits have been used in a lot of models to inquire about stability: a sufficient condition for stability of the process is the absorption of the fluid limits at 0. This is for instance used in \cite{rybsto} for two queues in tandem, and \cite{daifluid} for a multi-hop extension of Jackson networks. As stated in the next section, stability for QB-CSMA is already known but fluid limits can also be used to derive more detailed properties of the network. See for instance \cite{daifluidelay} where the authors derive the limits moments of the size of queues. With a method developed in \cite{brams98} and \cite{will98} based on the fluid properties of the model, Stolyar proved a Heavy traffic state space collapse for Max-Weights on a switch in \cite{Sto04} with a resource pooling assumption and later \cite{shahmw} removing this assumption. The goal of this paper is to prepare future work to study Queue-Based CSMA in a heavy traffic regime. 
	
	One of the weakness of the homogenization result used is that the bound obtained can be non-negligible for a fluid approximation when some queue lengths are too small. The complexity of QB-CSMA that was alleviated by the homogenization reappears when this occurs. For this reason, we are only able to prove convergence to a fluid limit in general when the limiting process does not approach zero in any coordinate. The reflection issue for a complete interference graph and the method associated to it are the main results of this paper. The method tailored to QB-CSMA can hopefully lead to more result in this vein.
	
	In a complete interference graph, we are able to describe accurately how the process is reflected. In short, if the initial value is not 0 for every queue, they all become positive instantaneously in the limit. This is due to the fact that no queue can be decreasing if it is too small compared to the sum of all queues and is done by a detailed analysis of pre-limit service times. From that point, either all queue reach zero at the same time if the arrival rates are sub-critical, or no queue reaches zero in finite time if the arrival rates are super-critical. The behaviour of the limit is given by an ODE. Near 0, the convergence is handled using tightness and the fact that the limit must be continuous. This reflection/absorption phenomenon is the most challenging and technical part of the article. This difficulty is due to the fact that QB-CSMA can idle even when no queue is empty. We need an had-hoc construction to prove this result. This reflection is in contrast with other functional limit theorem results where reflections are given by a Skorokhod map (introduced in \cite{sko1}, \cite{sko2} and generalized to multidimensional problems in \cite{tanaka}), see for instance \cite{daifluidelay}, \cite{brams98}, \cite{will98}, \cite{kangwill} for examples in queueing.

	\subsection{\large{\underline{Organization of the paper}}}

	The next section is dedicated to giving some notations and describing formally the model. Then, in Section \ref{sec:limitbrief} we begin by a brief description of the candidate fluid limits in both cases of interest (general interference graph with positive initial conditions and complete interference graph with general initial conditions). Proofs of existence and uniqueness of solution to the equations presented in this section are purely deterministic and can be found in Appendix \ref{app:ode}. Next, the main results are presented in Section \ref{sec:mainflui}. Some preliminary tightness and martingale results about the process of queue length are stated and proved in Section \ref{sec:preli}. We state the modified version of the homogenization result from \cite{us} in Section \ref{sec:loca} (the result can be found in \cite{thesis}). The proof of Theorem \ref{thm:main} is a direct consequence and is located in Section \ref{sec:thm1}. For the proof of Theorem \ref{thm:beyond}, some more work about reflections near the origin is required and the arguments can be found in Section \ref{sec:thm2}. 
	\subsection{\large{\underline{Notations and model description}}}\label{sec:model}
	
	We consider a finite set $V$ of $n<+\infty$ nodes. Each node $v \in V$ is an $M/M/1$ queue with the FCFS service discipline and vacation, jobs enter the queue along a Poisson process of intensity $\lambda_v > 0$. We denote by $Q_v(t) \in \mathds{N} \coloneqq \{0, 1, \ldots, \}$ the length of $v$'s backlog at time $t$ and by $(\sigma_v(t))_{v\in V} \in \{0,1\}^V$ the activity process: the server at $v$ is active and processing pending requests at unit rate whenever $\sigma_v(t) = 1$, and $\sigma_v(t) = 0$ otherwise. Put differently, $\sigma_v(t)$ is the instantaneous service rate of node $v$ at time $t$. We will later have to refer to an empty schedule (no queue active) and write this activity state $\textbf{0}$. Next, we define $\lambda \coloneqq (\lambda_v, v \in V)$, $Q(t) \coloneqq (Q_v(t), v \in V)$ and $\sigma(t) \coloneqq (\sigma_v(t), v \in V)$. For any $m<+\infty$, $x\in \mathds{R}^m$, we define $s(x)=\sum_{k=1}^m x_k$.
	
	We also introduce here a family of stopping times : for any $f:\reels_+\to \reels^V$ (we keep the same notation for real valued functions), and $\epsilon>0$,
	\[\tau^\epsilon(f)= \inf\lbrace t>0, \,  \min_v f_v(t)\leqslant \epsilon\rbrace.\]

	The nodes are placed on a simple undirected graph $G=(V,E)$. An edge between two nodes indicates that they cannot be active at the same time. This is used to model interference constraints in a wireless network. In other words, admissible schedules are \textit{stable sets} of $G$. The $\sim$ sign will be used to signify the existence of an edge ($v\sim w \Leftrightarrow \lbrace v,w\rbrace\in E$).  We will have two equivalent representations for the schedule: we will sometimes see it as a subset of nodes $\sigma\subset V$ when speaking of ``adding'' or ``removing'' a node in the schedule and otherwise as a vector of $\lbrace 0,1\rbrace^V$ by identifying nodes currently active in the schedule with non-zero entries of the service rate vector. A stable set of $G$ is given by $\sigma\in \lbrace 0,1\rbrace^V$ such that $v\sim w\Rightarrow\sigma_v+\sigma_w\leqslant 1$. The admissible service decisions are elements of 
	\[S(G)\coloneqq \lbrace \sigma\in \lbrace0,1\rbrace^V \mid v\sim w\Rightarrow\sigma_v+\sigma_w\leqslant 1 \rbrace .\] 
	Given the current schedule $\sigma$, the queue-length process $Q$ evolves as $n$ independent $M/M/1$ queues with service rates $\sigma$, input rates $\lambda$ and FCFS service discipline. On the other hand, $\sigma$ also evolves: given the queue-length process $Q$, an active node $v$ with $\sigma_v = 1$ deactivates at rate $\Psi_-(Q_v)$ for some deactivation function $\Psi_-$, and an inactive node $v$ with $\sigma_v = 0$ activates at rate $\Psi_+(Q_v)$ for some activation function $\Psi_+$, provided no neighbouring node is active.
	
	To be more formal, $(Q, \sigma)$ is a Markov process on $\mathds{N}^V \times \{0,1\}^V$ with infinitesimal generator $L$ that can be decomposed as the sum of two generators:\begin{itemize}
		\item the generator $L^\sigma_{\mathrm{s}}$ of the \textit{slow} queue-length process $Q$ whose dynamic depends on $\sigma$;
		\item and the generator $L^q_{\mathrm{f}}$ of the \textit{fast} activity process $\sigma$ whose dynamic depends on $q$.
	\end{itemize}
	The terminology \textit{slow} and \textit{fast} refers to the time scales in the  stochastic averaging principle from Section \ref{sec:sap}. Thus, $L$ acts on functions $f: \mathds{N}^V \times S(G) \to \mathds{R}$ as
	\[ L[f](\sigma,q) = L_{\mathrm{s}}^\sigma[f(\sigma, \cdot)](q) + L_{\mathrm{f}}^q[f(\cdot,q)](\sigma) \]
	with
	\begin{equation} \label{eq:slow-L}
	L_{\mathrm{s}}^\sigma[g](q) = \sum_{v \in V} \lambda_v \left( g(q+e^v) - g(q) \right) + \sum_{v \in V} \sigma_v \ind_{q_v>0} \left( g(q-e^v) - g(q) \right)
	\end{equation}
	and
	\begin{multline} \label{eq:fast-L}
	L_{\mathrm{f}}^q[h](\sigma) = \sum_{v \in V} \sigma_v \Psi_-(q_v) \left( h(\sigma - e^v) - h(\sigma) \right)\\
	+ \sum_{v \in V}\prod_{w \sim v} (1-\sigma_w)(1-\sigma_v)  \Psi_+(q_v) \left( h(\sigma + e^v) - h(\sigma) \right)
	\end{multline}
	with $g: \mathds{N}^V \to \mathds{R}$ and $h: S(G) \to \mathds{R}$ arbitrary functions and $e^v \in \{0,1\}^V$ with $0$'s everywhere except at the $v$th coordinate equal to $1$. One can check that for any $q\in \entier^V$, $L^q_{\mathrm{f}}$ admits a unique reversible distribution denoted $\pi^q$. Because of the homogenization result of \cite{us}, we consider polynomial activation rates
	\[ \Psi_+(x) = \frac{(x+1)^a}{1 + (x+1)^a} \in [0,1] \ \text{ and } \ \Psi_-(x) = 1 - \Psi_+(x), \ x \in \mathds{N}, \]
	with $a > 0$ the parameter of this algorithm. Polynomial rates are the fastest growing functions such that homogenization is possible. In this case, $\pi^q$ is given by
	
	\[ \pi^q(\sigma) = \frac{\prod\limits_{v\in V}(1+q_v)^{a\sigma_v}}{\sum\limits_{\rho\in S(G)} \prod\limits_{v\in V}(1+q_v)^{a\rho_v}},  \sigma \in S(G). \]
	
	Let $\ell^q$ be the spectral gap of $L_{\mathrm{f}}^q$ and assume that $0<\beta<\infty$ is such that
	\[\ell^q\geqslant C\norm{q+1}{\infty}^{-a\beta}.\]
	The existence of such a $\beta$ is guaranteed for instance by \cite{shashin}, Lemma 3. Their result states that $\beta\leqslant 2(n+1)$ but a better bound can be obtained for specific topologies. 
	
	An important tool to understand the behaviour of the network is the so-called ``homogenized process'' where the service rates are replaced by their steady state distribution. Let $g:\entier^V\to \reels$. The generator is given by 
	\begin{equation}
	L_{\mathrm{h}} [g](q)=  \sum_{v \in V} \lambda_v \left( g(q+e^v) - g(q) \right) + \sum_{v \in V} \pi^{q}(\sigma_v=1) \ind_{q_v>0} \left( g(q-e^v) - g(q) \right).\label{eq:homo-L}
	\end{equation}
	
	When all queue lengths are large, only the stable sets of largest size will matter in the invariant measure.  For any $\sigma\in S(G)$, $\left\lvert \sigma \right\rvert$ is its size (the number of active nodes). We can define the constant $\Upsilon\coloneqq \max_{\sigma \in S(G)} \lvert \sigma \rvert$ (we omit the dependency in $G$) and define
	
	\[S^*\coloneqq \left\lbrace \sigma\in S(G) \mid \left\lvert \sigma \right\rvert=\Upsilon\right\rbrace .\]
	
	When $q_v$ is larger than $0$ for every $v$, for a large parameter $N$, $\pi^{Nq}\approx \pi_\infty^q$ given by 
	
	\[\pi_\infty^q(\sigma)=0 \text{ if } \sigma\notin S^* \text{ and }\pi_\infty^q(\sigma)= \frac{\prod\limits_{v\in V}q_v^{a\sigma_v}}{\sum\limits_{\rho\in S^*} \prod\limits_{v\in V}q_v^{a\rho_v}} \text{ if } \sigma\in S^*.\] The instantaneous service rate  will be given by $\overline{\pi}^q(v)=\pi_\infty^q(\sigma_v=1)$ for all $v\in V$. There is a partially uniform result:
	\begin{lemma}\label{lem:servicerate}
		
		Let $C_->0$ and $C_+<+\infty$. We have 
		\begin{equation}
		\forall v\in V,\; \sup\limits_{C_-\leqslant \min_w q_w\leqslant \max_w q_w\leqslant C_+, 
			q\in \frac{1}{N}\entier^V} \lvert\pi^{Nq}(\sigma_v=1) - \overline{\pi}^q(v)\rvert \to 0 \text{ as } N\to +\infty. \label{eq:servicerate}
		\end{equation}
	\end{lemma}
	
\textbf{	{\underline{Proof sketch:}}}
		The full proof is omitted to lighten the paper. The general idea is to decompose $\pi^{Nq}(\sigma_v=1)$ along the size of independent sets. Because of the scaling, independent set of size smaller than $\Upsilon$ will be negligible in the limit. A non uniform result can be obtained by noticing that for any $q$ bounded and bounded away from zero, $\pi^{Nq}(\sigma_v=1)\to \bar{\pi}(v)$. In order to obtain a uniform result, we need to use a Taylor expansion, which will provide explicit bounds depending on $C_-,C_+$ and $N$ that vanish as $N\to +\infty$.
	\qed

	\section{\large{Main results and heuristic}}
	\subsection{\large{\underline{Limiting process}}}\label{sec:limitbrief}
	We define here briefly the limiting processes that we will encounter for the fluid limits. More details about existence and uniqueness are given in Section \ref{sub:limg}, proofs of the results are given in Appendix \ref{app:ode}.
	
	Let
	\begin{equation}\label{eq:g}
	\begin{array}{ll}
	g:\,\reels_+^V\setminus\lbrace \textbf{0}\rbrace &\to [-1,1]^V\\
	\quad \quad q &\longmapsto \lambda-\overline{\pi}^q
	\end{array}
	\end{equation}
	be the difference between arrival rates and homogenized departure rates.  Let $q^0\in (0,+\infty)^V$. We will consider the equation with unknown $f$ with values in $(0,+\infty)^V$ and domain included in $[0,+\infty)$: \begin{equation}
	\left\lbrace\begin{array}{ll}
	f'=g(f) & \\
	f(0)=q^0 \end{array}\right. . \label{eq:limit}
	\end{equation}
	\begin{lemma}\label{lem:uniquegen}
		For any $q^0\in (0,+\infty)^V$, there exists at least one solution of \eqref{eq:limit} defined in an open interval $D\subset\reels$ containing $0$.\\
		If $D_1$ and $D_2$ are two open intervals containing 0 such that $q^1:D_1\to (0,+\infty)^V$ and $q_2:D_2\to (0,+\infty)^V$ are two solutions of \eqref{eq:limit}, then $q^1(t)=q^2(t)$ for any $t\in D_1\cap D_2$.

	\end{lemma}
\textbf{{\underline{Proof:}}}

		Since $g$ is locally Lipschitz on $(0,+\infty)^V$, this is a direct application of the Cauchy Lipschitz Theorem, see Chapter 5 Section 3.1 and 3.3 of \cite{ode}.  
	\qed 
	
	We call $q^*(\cdot, q^0)$ the unique maximal solution to this ODE. This type of uniqueness is standard in ODE (\cite{ode} Chapter 5). It is known that the solution is unique up to the time it exits $(0,+\infty)^V$ so for any initial condition, so there exists a maximal open set $D^*(q^0)$ containing 0 on which to define the solution. By standard results, $\sup D^*(q^0)=+\infty$ or the solution exits $(0,+\infty)^V$ in finite time. The exit time only depends on the initial condition and we denote it 
	\[T_{\text{ext}}(q^0)\coloneqq \sup D^*(q^0).\] 
	\begin{definition}\label{def:ext}
		For any $q^0\in (0,+\infty)^V$, we call
		\[ q^*(\cdot,q^0):[0,+\infty)\to \reels_+^V\] 
		the restriction on $[0,+\infty)$ of the maximal solution of \eqref{eq:limit} with initial condition $q^0$. The solution is extended after $T_{\text{ext}}(q^0)$ by stating that 
		\[\forall s\geqslant 0,\; q^*(T_{\text{\textnormal{ext}}}(q^0)+s,q^0)=\lim_{t\to T_{\text{\textnormal{ext}}}(q^0)}q^*(t,q^0) .\]
	\end{definition}
	For a complete interference graph, let's define the initial value problem
	\begin{equation}
	\left\lbrace\begin{array}{ll}
	f'=g(f)   &\text{if } f\neq\textbf{0} \\
	s(f)(t)=\max(s(q^0)+(s(\lambda)-1)t,0)\\
	f(0)=q^0 \end{array}\right. \label{eq:limitCIG}
	\end{equation}
	and denote its solutions $q^*_c(\cdot,q^0)$.

	\subsection{\large{\underline{Main result}}}\label{sec:mainflui}
	The first theorem of this paper concerns the fluid limit of the system. It states that over small enough time horizon, the scaled process of queue lengths converges uniformly in probability to $q^*$ given in the previous section. Recall the definition of the process 
	\[Q^N= \frac{Q(N\cdot)}{N},\]
	and $\beta$ is such that 
	\[\ell^q\geqslant C\norm{q+1}{\infty}^{-a\beta}.\]
	For instance $\beta=1$ in the case of a complete interference graph, see \cite{us} Lemma 5.3.
	
	\begin{theorem}\label{thm:main}
		Assume that the two following assumptions hold:
		\begin{itemize}
			\item $2a\beta < 1$ ;
			\item $Q^N(0) \to q^0$ for some $q^0\in (0,+\infty)^V$.
		\end{itemize}
		Then $Q^N(\cdot)\overset{\P}{\to} q^*(\cdot,q^0)$ as $N \to +\infty$ uniformly over compact sets of $[0,T_{\text{ext}}(q^0))$.
	\end{theorem}

	Because of the special structure of complete interference graph ( illustrated by Lemma \ref{lem:qpos}), this theorem is enough to state that in this case, if $s(\lambda)\geqslant 1$ and $Q^N(0)\to q^0$ with $q^0_v>0$ for every $v\in V$, \[Q^N\overset{\P}{\to} q^*(\cdot,q^0)\] uniformly over compact time intervals: in this case, no queue will ever reach zero. The comparison to an ODE is troublesome when some queues start empty or after they reach $0$ because homogenization may fail at that time. This does not happen in the super critical case. In general, for any interference constraints and initial condition such that a solution to the limiting ODE never exits $(0,+\infty)^V$, this theorem is enough to prove convergence of $Q^N$ uniformly over compact of $[0,+\infty)$.

	When $G$ is a complete graph, it is possible to control potential reflections at the boundary of $(0,+\infty)^V$ separately to obtain the following result:
	\begin{theorem}\label{thm:beyond}
		Let $q^0\in [0,+\infty)^V$ and $\lambda\in \reels_+^V$ be fixed. Assume the following: \begin{itemize}
			\item $Q^N(0)\to q^0$ as $N\to +\infty$,
			\item  $G$ is a complete interference graph,
			\item and $a<\frac{1}{2}$.
		\end{itemize} 
		Then $Q^N \overset{\P}{\to} q^*_c$ uniformly over compact time interval for $q^*_c$  characterized as the unique solution to \eqref{eq:limitCIG}.
	\end{theorem} 
	
	The proof of the first theorem relies heavily on two essential results stated in the next section: \ref{sec:preli}. First $Q^N$ is tight and any limiting point is a continuous function. There does not lie the difficulty and the proof is straightforward. Second  the homogenization result Corollary \ref{cor:scaling}, an application of the homogenization result of \cite{thesis}.
	
	Theorem \ref{thm:main} is a direct consequence of these results. Using homogenization, we are able to bound the distance between $Q^N$ and its homogenized version until $q^*$ reaches the boundary of its state space. Continuity of any limiting point and the fact that it cannot approach zero will lead us to the result. For a complete interference graph, we prove in Lemma \ref{lem:pos} that when starting with some null queues, all queues become positive for positive time and then use property of the limiting process to prove convergence.
	
	\subsection{\large{\underline{Preliminary results}}}
	We will begin this section by proving a tightness result for the process of queue lengths. Then, we will spend the next subsection defining the localization set and explaining the homogenization result of \cite{us,thesis}. Finally, we will state the existence and uniqueness result for \eqref{eq:limit} and \eqref{eq:limitCIG}. Recall that $C$ denotes a numerical constant allowed to depend on $a$, $n$, $\lambda$, and $G$. 
	\subsubsection{\large{\underline{Tightness of $Q^N$}}}\label{sec:preli}

	We begin with a classical result consequence of the martingale problem:
	\begin{lemma}	\label{lem:mart}
		For any $v\in V$, we define the process $M^N_v$ as
		\[\forall t\geqslant 0,\; M^N_v(t)=Q^N_v(t)-\int_{0}^{t}\left(\lambda_v-\sigma^N_v(s) \right)\d s.\]
		For any $v\in V$, $M^N_v$ is a martingale of bracket $\langle M_v\rangle_t=\int_{0}^t \dfrac{\lambda_v+\sigma_v^N(s)}{N}\d s$ and
		\begin{equation}
		\cesp{}{\sup\limits_{t\leqslant T}\lvert M^N_v(t)\rvert}\leqslant 2\sqrt{\frac{T(\max_v\lambda_v+1)}{N}}\to 0\text{ as }N\to +\infty.\label{eq:mart}	\end{equation} 
	\end{lemma}
	
\textbf{{\underline{Proof:}}}

		The fact that $M^N_v$ are martingales and the expression of their bracket are consequences of the martingale problem, see for instance \cite{ethierkurtz} chapter IV (Proposition 1.7) for more information in a more general setting. In addition,
		\begin{align*}
		\cesp{}{\sup\limits_{t\leqslant T}\lvert M^N_v(t)\rvert}&\overset{(a)}{\leqslant} \sqrt{\cesp{}{\sup\limits_{t\leqslant T} M^N_v(t)^2}}\\
		&\overset{(b)}{\leqslant} 2\sqrt{\cesp{}{M^N_v(T)^2}}\\
		&\overset{(c)}{=}2\sqrt{\cesp{}{N\int_{0}^{T}\frac{\lambda_v+\sigma^N_v(s)}{N^2}\texttt{d}s}}\leqslant 2\sqrt{\frac{T(\max_v\lambda_v+1)}{N}}.
		\end{align*}
		
		Inequality $(a)$ comes from Jensen's inequality, $(b)$ comes from Doob's inequality and $(c)$ comes from the definition of the bracket.
	\qed 
	
	\begin{proposition}
		For any tight sequence initial condition $Q^N(0)$, $Q^N$ is tight for the topology of uniform convergence over compact time sets. Any limiting point is a distribution on continuous functions.\label{prop: ctight}
	\end{proposition}
	
\textbf{{\underline{Proof:}}}

		To obtain this property, recall the definition of the modulus of continuity of $Q^N$ from Chapter 2 of \cite{bill}: for any $\delta>0$, $T<+\infty$
		\begin{align*}
		\omega_{Q^N}(\delta)&=\sup\limits_{s,t\leqslant T, \lvert t-s\rvert\leqslant  \delta, v\in V} \lvert Q^N_v(t)-Q^N_v(s) \rvert\\
		&=\sup_{s,t\leqslant T, \lvert t-s\rvert \leqslant \delta, v\in V} \lvert\int_{s}^{t}(\lambda_v-\sigma^N_v(u)) \d u+M^N_v(t)-M^N_v(s)\rvert\\
		&\leqslant \delta\max_v (\lambda_v\vee \lvert\lambda_v-1\rvert)+2\sup\limits_{t\leqslant T, v\in V} \lvert M^N_v(t)\rvert
		\end{align*} 
		To conclude the proof of this proposition, use Markov inequality on ~\eqref{eq:mart} to obtain: 
		\begin{equation}
		\cpro{}{	\omega_{Q^N}(\delta)\geqslant \delta(\max_v (\lambda_v\vee \lvert\lambda_v-1\rvert)+\nu)}\to 0 \text{ as } N\to +\infty \text{ for any } \nu>0,\label{eq: contmod}
		\end{equation}
		which coupled with convergence of the initial condition implies the result with Theorem 7.3 in \cite{bill}.
	\qed 
	
	\subsubsection{\large{\underline{Localization and homogenization}}}\label{sec:loca}

	Let us define a localization set of the form $U \subset \reels^V_+$ such that
	\begin{equation}
	U \coloneqq \left \lbrace q \in \reels_+^V: \forall v \in V,\, q_v \in (C_-,C_+) \right \rbrace,\label{eq:uloc}
	\end{equation}
	Recall $q^*$ the solution of \eqref{eq:limit} or \eqref{eq:limitCIG} depending on the context, and 
	\[\tau^\epsilon(f)=\inf\left\lbrace t>0,\, \min_{v\in V} f_v(t)\leqslant \epsilon\right\rbrace.\] By continuity of $q^*$, for any $T<+\infty$, there exists a $C_T$ such that
	
	\[\sup_{t\leqslant T}\max_v q_v^*(t)<C_T.\]
	
	For any fixed $q^0\in (0,+\infty)^V$ and $0<\epsilon<\min_v q^0_v$, let us fix $T$ a finite horizon such that $T\leqslant \tau^\epsilon(q^*)$. To obtain convergence uniformly over compact sets contained in $[0,\tau^\epsilon(q^*))$ define $U^\epsilon$ in \eqref{eq:uloc} with $C_-=\frac{\epsilon}{3}$ and $C_+=3C_{\tau^\epsilon(q^*)}$ (we emphasize the dependency in $\epsilon$ only here but omit it until needed in the rest of the paper).
	Introduce the exit time of $Q^N$ from $U$:
	\[ \bar{\tau}^{U}(Q^N)=\inf \left \lbrace t \geq 0: Q^N(t) \notin U \right \rbrace. \]

	Define also 
	
	\[ T_{\text{\textnormal{tube}}}^N\coloneqq  \inf\left\lbrace t>0 , \; \norm{Q^N(t)-q^*(t)}{\infty}>\min(C_-,C_+/3) \right\rbrace.\]
	
	We present here a lemma used to remove the localization by replacing it with a ``tube condition'' around the limiting process. The behaviour of the limiting process allows us to remove this tube condition as well in some cases.
	\begin{lemma}\label{lem:m-M}
		Almost surely $T_{\text{\textnormal{tube}}}^N < \bar{\tau}^{U}(Q^N)$. In particular, $Q^N(t \wedge T_{\text{\textnormal{tube}}}^N) \in U$ for all $t\leqslant T$.
	\end{lemma}
\textbf{{\underline{Proof:}}}
		Recall that by definition of the finite time horizon 
		\[\epsilon<\inf_{0\leqslant t\leqslant T, v\in V} q^*_v(t) \leqslant \sup_{0 \leqslant t \leqslant T, v\in V }q^*_v(t)< C_{\tau^\epsilon(q^*)},\]
		or equivalently 
		\[3C_-<\inf_{0\leqslant t\leqslant T, v\in V} q^*_v(t) < \sup_{0 \leqslant t \leqslant T, v\in V }q^*_v(t)\leqslant C_+/3.\]
		Since the jump size of $Q^N$ is $\frac{1}{N}$ we have $\sup\limits_{t\leqslant T_{\text{\textnormal{tube}}}^N}\norm{Q^N(t)-q^*(t)}{\infty}\leqslant \min (C_-, C_+/3)+\frac{1}{N}$. In addition, for large $N$, $\frac{1}{N}$ is smaller than $\min (C_-, C_+/3)$. A direct consequence is that
		\[3C_--2C_-< \inf_{0\leqslant t\leqslant T_{\text{\textnormal{tube}}}^N\wedge T}Q^N_v(t)\leqslant  \sup_{0\leqslant t\leqslant T_{\text{\textnormal{tube}}}^N\wedge T}Q^N_v(t)< 3C_+/3.\]
		Said otherwise, for any $T\leqslant \tau^\epsilon(q^*)$, 
		\[\left\lbrace \bar{\tau}^U(Q^N)\geqslant T\right\rbrace \subset \left\lbrace T_{\text{\textnormal{tube}}}^N<\bar{\tau}^{U}(Q^N),\, \bar{\tau}^U(Q^N)\geqslant T\right\rbrace.\]
	\qed 
	
	Using the localization set from \eqref{eq:uloc}, we are able to state the homogenization result that is used in this paper. The method is derived from \cite{us} Proposition 3.4, or see \cite{thesis}, Corollary 3.7 for a proof of this result.

	\begin{corollary} \label{cor:scaling}
		Assume that 
		\[\ell^q\geqslant C\norm{q+1}{\infty}^{-a\beta}.\] 
		For any $v \in V$, $t\leqslant T$, if $a\beta<1/2$ we have
		\begin{equation}
		\E \left[ \sup_{0 \leqslant t \leqslant T \wedge \bar{\tau}^U(Q^N)} \left \lvert \int_0^t \left( \sigma^N_v(s)-\pi^{NQ^N(s)}(v) \right)\d s \right \rvert \right]\leqslant CN^{a\beta-1/2}\log(N)^{3/2} .\label{eq:homogenization}
		\end{equation}
	\end{corollary}
	Without the condition $a\beta<1/2$, the bound is a little bit different but the range of $a$ such that it is negligible does not change. Interestingly, the method from \cite{Kurtz92} would not allows for such a result but \cite{Luc13} (Theorem 1.5) would yield a similar bound. The reason behind this new homogenization bound is given in \cite{us} and will hopefully lead to new heavy traffic results.
	
\textbf{\underline{Proof sketch:}}
		The main idea is to write $\sigma_v-\pi^{q}(v)$ in terms of $\phi_v(q,\sigma)$, the solution to the Poisson equation :
		\begin{equation} \label{eq:Poisson}
		L^{q}_{\mathrm{f}} [\phi_v](\sigma) = \sigma_v - \pi^q[v],\; \pi^q[\phi_v]=0,
		\end{equation}
		and control the term we obtain using the martingale problem for $\phi_v(Q^N,\sigma^N)$. We have to use a generalization of Lemma \ref{lem:mart} with $\phi_v(q,\sigma)$ instead of $q_v$ but this is covered by \cite{ethierkurtz} Proposition 1.7.  More formally,
		\begin{align*}
		\int_0^{Nt} \left( \sigma_v(s)-\pi^{Q(s)}(v) \right)\d s&=\int_{0}^{Nt}	L^{Q(s)}_{\mathrm{f}}[\phi_v(\cdot,Q(s))](\sigma(s))ds\\
		&=\phi_v(Q(Nt),\sigma(Nt))-\phi_v(Q(0),\sigma(0))-\int_{0}^{Nt}L^{\sigma(s)}_{\mathrm{s}}[\phi_v(\sigma(s),\cdot)](Q(s))ds\\
		&\quad+M_v(Nt),
		\end{align*}
		where $M_v$ is a martingale of bracket given by
		\begin{equation*}
		\langle M_v\rangle_t=\int_{0}^{t}\Gamma[\phi_v](Q(s),\sigma(s)),
		\end{equation*}
		with \begin{align*} \label{eq:Gamma}
		\Gamma [f](q,\sigma) &= \sum_{v \in V} \lambda_v \left( f \left( q+e^v,\sigma \right)-f(q,\sigma)\right)^2\\
		& \hspace{-5mm} +  \sum_{v \in V} \sigma_v\ind_{q_v>0} \left(f \left( q-e^v,\sigma \right)-f(q,\sigma)\right)^2 \notag\\
		& \hspace{-5mm} +  \sum_{v \in V} \sigma_v \Psi_-(q_v) \left( f(q,\sigma - e^v) - f(q,\sigma) \right)^2\\
		& \hspace{-5mm} + \sum_{v \in V}\prod_{w \sim v} (1-\sigma_w)(1-\sigma_v)  \Psi_+(q_v) \left( f(q,\sigma + e^v) - f(q,\sigma) \right)^2.
		\end{align*} Explicit bounds on the regularity of $\phi_v$ as a function of $q$ are needed to control the bracket. The result then follows from the assumption on the spectral gap and Lemma 5.2 of \cite{us}.
	\qed 
	\subsubsection{\large{\underline{Existence and uniqueness of solutions to \eqref{eq:limit} and \eqref{eq:limitCIG}}}}\label{sub:limg}
	
	Recall the definition of $S^*$:
	\[S^*= \left\lbrace \sigma\in S(G) \mid \left\lvert \sigma \right\rvert=\Upsilon\right\rbrace,\] 
	and the definition of the stopping times $\tau^\epsilon$:
	\[\tau^\epsilon(f)= \inf\lbrace t>0, \,  \min_v f_v(t)\leqslant \epsilon\rbrace.\]
	We begin by identifying a potential fluid limit using the homogenized process. On the homogenized process, because of Lemma \ref{lem:servicerate}, we will localize the process such that the asymptotic service rate will be given for any $v\in V$ by
	\[ \bar{\pi}^q(v)= \pi_\infty^q(\sigma_v=1)= \sum_{\sigma\in S^*}\frac{\sigma_v\prod\limits_{w\in V}q_w^{a\sigma_w}}{\sum\limits_{\rho\in S^*} \prod\limits_{w\in V}q_w^{a\rho_w}} ,\]
	recall as well the definition of $g$ from \eqref{eq:g}:
	\begin{equation*}
	\begin{array}{ll}
	g:\,\reels_+^V\setminus\lbrace \textbf{0}\rbrace &\to [-1,1]^V\\
	\quad \quad q &\longmapsto \lambda-\overline{\pi}^q
	\end{array}
	\end{equation*} 
	
	The function $g$ is locally Lipschitz on $(0,+\infty)^V$ by differentiability.

	For a complete interference graph, we would like to be able to prove convergence to a unique fluid limit for any initial condition, possibly null even though $g(\textbf{0})$ is not defined. Starting from \textbf{0} does make sense for the ODE given in the previous section but the complete interference graph allows us to make statements. We introduce the next equation in order to specify what happens when approaching \textbf{0} or starting in this state. When $G$ is a complete interference graph, let's define $q_c^*$ as the only (existence and uniqueness are proved in Lemma \ref{lem:uniquecig}) function $f:\reels_+\to \reels^V_+$ that is continuous, differentiable almost everywhere (everywhere except at $t=0$ if $q^0=\textbf{0}$ or $\inf\left\lbrace t>0, f(t)=\textbf{0}\right\rbrace$), such that 
	\begin{equation*}
	\left\lbrace\begin{array}{ll}
	f'=g(f)   &\text{if } f\neq\textbf{0} \\
	s(f)(t)=\max(s(q^0)+(s(\lambda)-1)t,0)\\
	f(0)=q^0 \end{array}\right.
	\end{equation*}
	We now give an informal description of two important lemmas that will justify this definition and link solutions to this initial value problem with solutions to \eqref{eq:limit}. The lemmas are proved in Appendix \ref{app:ode}
	
	\begin{lemma}\label{lem:uniquecig}
		Let $q^N\in (0,+\infty)^V$ such that $q^N\to q^0\in \reels_+^V$. 
		\begin{itemize}
			\item If $q^0_v>0$ for all $v\in V$, $q^*(\cdot,q^N)\to q^*(\cdot,q^0)$ uniformly on compact time sets .
			\item  If there exists $v\in V$ such that $q^0_v=0$, the limit is still well defined and \[q_c^*(\cdot,q^0)\coloneqq\lim_{N\to +\infty}q^*(\cdot,q^N) \] is a solution to the initial value problem \eqref{eq:limitCIG}.
			\item For any $q^0\in \reels^V_+$, the solution to the initial value problem \eqref{eq:limitCIG} is unique.
		\end{itemize}
	\end{lemma}
	In other words, the additional information about the sum of coordinate is enough to uniquely specify a solution to the initial value problem \eqref{eq:limitCIG} among the solutions to the ODE \eqref{eq:limit} even when the solution may not be unique. In addition, even if $q^0=\textbf{0}$ we can construct a solution to \eqref{eq:limitCIG} using solutions to \eqref{eq:limit} for positive initial values $q^N\to \textbf{0}$. For this type of solution we can also define an exit time as 
	\[T_{\text{ext}}(q^0)\coloneqq\lim_{N\to +\infty}T_{\text{ext}}(q^N).\]

	Uniqueness in this lemma is in the same sense as in Lemma \ref{lem:uniquegen}. To prove this lemma, we explain in the next lemma why the complete interference graph constitutes a particular case: the dynamic \eqref{eq:limitCIG} ensures that all coordinates touch $0$ at the same time if they do at all. Recall the definition of the stopping times  
	\[\tau^0(f)=\inf\lbrace t>0, \; \min_v f_v(t)\leqslant 0\rbrace.\]
	The next lemma is required to prove uniqueness of solutions of \eqref{eq:limitCIG}. We will justify the existence of solutions in Lemma \ref{lem:uniquecig} but it is in fact guaranteed by the Cauchy Peano Arzelà theorem (see\cite{ode}). We will prove the next lemma for any solution to \eqref{eq:limitCIG} without assuming uniqueness.
	\begin{lemma}\label{lem:qpos}
		Let $q^0\in \reels_+^V$ and $q_c^*$ any solution to \eqref{eq:limitCIG}. Then \label{lem:simul} $\tau^0(q_c^*)=\tau^0(s\circ q_c^*)$. In addition,
		\begin{itemize}
			\item If $s(\lambda)> 1$ , $\tau^0(s\circ q^*_c)=+\infty$. 
			\item If $s(\lambda)<1$, $\tau^0(s\circ q^*_c)=\frac{s(q^0)}{1-s(\lambda)}$
			\item If  $s(\lambda)=1$ and $q^0\neq \text{\textbf{0}}$, $\tau^0(s\circ q^*_c)=+\infty$.
			\item  If $s(\lambda)=1$ and $q^0= \text{\textbf{0}}$, $\tau^0(s\circ q^*_c)=0$.
		\end{itemize}
	\end{lemma}
\textbf{{\underline{Proof sketch:}}}
		The general idea is that no coordinate can be decreasing when too small compared to the sum. The expressions for $\tau^0(s\circ q^*_c)$ are direct consequences of the description of the sum, which has an easy to use expression. See Appendix \ref{app:ode} for a proof of this result.
	\qed 
	Because of this Lemma, we keep the $q^*(\cdot,q^0)$ notation for solutions to \eqref{eq:limitCIG}  when unambiguous: they coincide with the extension given in Definition \ref{def:ext} when $q^0\in (0,+\infty)^V$ and $q^*_c$ is the only solution of \eqref{eq:limit}  defined on $[0,+\infty)$ when $q^0_v=0$ for some $v$.

	\section{\large{General interference graph up to $\tau^0(q^*)$, Theorem \ref{thm:main}}} \label{sec:thm1}
	
	Fix $\epsilon>0$ and $T\leqslant \tau^\epsilon(q^*)$. Recall the definitions:
	\[U \coloneqq \left \lbrace q \in \reels_+^V: \forall v \in V,\, q_v \in (C_-,C_+) \right \rbrace,\]
	with 
	\[ (C_-,C_+)=(\dfrac{\epsilon}{3}, \max_{t\leqslant \tau^{\epsilon}(q^*), v\in V}q^*_v(t))\]
	and $C_->0$ and $C_+<+\infty$. To prove Theorem \ref{thm:main}, we will establish its equivalent for the stopped process $Q^N( \, \cdot \,\wedge T_{\text{\textnormal{tube}}}^N\wedge T)$ using Gronwall's lemma. We then transfer the result on the stopped process to $Q^N(\cdot\wedge T)$ using Lemma~\ref{lem:m-M}. This gives us uniform convergence over compact sets of $[0,\tau^\epsilon(q^*))$ for any $\epsilon>0$. We also give the arguments to extend up to $\tau^0(q^*)$ at the end of the proof. 
	
\textbf{{\underline{Proof of Theorem \ref{thm:main}:}}}

		Recall that for any $t \leqslant T\leqslant \tau^\epsilon(q^*)$,
		\[q^*(t)=q^0+\int_{0}^{t}\left(\lambda-\overline{\pi}^{q^*(s)}\right)\texttt{d}s.\]
		Using Lemma \ref{lem:mart}, we get
		\[Q^N(t)=Q^N(0)+\int_{0}^{t}\lambda-\ind_{Q^N(s)>0}\sigma^N(s)\d s +M^N(t).\]
		Adding and subtracting \[\int_{0}^{t}\left(\pi^{NQ^N(s)}(\sigma_v=1)+\bar{\pi}^{Q^N(s)}(v)\right)\d s,\] we get
		\begin{align*}
		Q^N_v(t)-q_v^*(t)&=Q^N_v(0)-q^0_v+\int_{0}^{t}\left(\overline{\pi}^{q^*(s)}(v)-\overline{\pi}^{Q^N(s)}(v)\right)\texttt{d}s+M^N_v(t) \\
		&+\int_{0}^{t}\left(\overline{\pi}^{Q^N(s)}(v)-\pi^{NQ^N(s)}(\sigma_v=1)+\pi^{NQ^N(s)}(\sigma_v=1)-\ind_{Q^N_v(s)>0}\sigma_v^N(s)\right)\texttt{d}s.\notag
		\end{align*}
		Define 
		\[\Theta_v^N= \left\lvert Q^N_v(0)-q^0_v\right\rvert+\eta_v^N(T)+h_v^N(T)+\mu^N_v(T),\]
		were 
		\[ \eta_v^N(T)=\sup\limits_{t\leqslant T}\left\lvert\int_{0}^{t\wedge T_{\text{\textnormal{tube}}}^N}\left(\overline{\pi}^{Q^N(s)}(v)-\pi^{NQ^N(s)}(\sigma_v=1)\right)\texttt{d}s\right\rvert,\]
		
		\[h^N_v(T)=\sup\limits_{t\leqslant T}\left\lvert\int_{0}^{t\wedge T_{\text{\textnormal{tube}}}^N}\left((\pi^{NQ^N(s)}(\sigma_v=1)-\ind_{Q^N_v(s)>0}\sigma_v^N(s)\right)\texttt{d}s\right\rvert\]
		
		and 
		
		\[\mu^N_v(T)=\sup\limits_{t\leqslant T\wedge T_{\text{\textnormal{tube}}}^N}\left\lvert M^N_v(t)\right\rvert.\]
		
		We get for all  $v\in V$,
		
		\[	\left\lvert Q^N_v(t\wedge T_{\text{\textnormal{tube}}}^N)-q_v^*(t\wedge T_{\text{\textnormal{tube}}}^N)\right\rvert\leqslant \Theta^N_v+\int_{0}^{t\wedge T_{\text{\textnormal{tube}}}^N}\left\lvert\overline{\pi}^{q^*(s)}(v)-\overline{\pi}^{Q^N(s)}(v)\right\rvert \texttt{d}s  .\]
		
		Remark $f: U\to \reels^V$ defined by $f(q)=\overline{\pi}^q$ is Lipschitz. So there exists $C<+\infty$ such that for any $q,q'\in U$
		\[\norm{f(q)-f(q')}{} \leqslant C \norm{q-q'}{},\]
		and thus,
		\[\int_{0}^{t}\left\lvert\overline{\pi}^{q^*(s)}(v)-\overline{\pi}^{Q^N(s)}(v)\right\rvert \texttt{d}s\leqslant C\int_{0}^{t}\norm{Q^N(s)-q^*(s)}{}\texttt{d}s.\]
		To sum up, for any $t\leqslant T\wedge T_{\text{\textnormal{tube}}}^N$, by summing over $v$, and denoting $\Theta^N=\sum_{v \in V} \Theta^N_v$,
		\[	\norm{ Q^N(t)-q^*(t)}{1}\leqslant \Theta^N+C\int_{0}^{t}\norm{Q^N(s)-q^*(s)}{1}\texttt{d}s . \]
		
		By Gronwall lemma (\cite{ethierkurtz}, Appendix 5), for any $t\leqslant T\wedge T_{\text{\textnormal{tube}}}^N$, we get \[\norm{ Q^N(t)-q^*(t)}{1}\leqslant\Theta^N\exp(Ct)\leqslant \Theta^N\exp(CT).\]
		By Lemma \ref{lem:servicerate} for all $v\in V$, $\cesp{}{\eta^N_v(T)}\to 0$ as $N\to +\infty$. Further, because of Lemma \ref{lem:m-M}, $T^N_{tube}<\bar{\tau}^U(Q^N)$ so before $T^N_{tube}$ all queue lengths are bounded away from zero. Corollary \ref{cor:scaling} states that $\cesp{}{h^N_v(T)}\to 0$ as $N\to +\infty$ if $a\beta<1/2$. Every $\mu^N_v$ converges to $0$ in the mean, by Lemma \eqref{lem:mart}. Finally, $Q^N_v(0)-q^0_v\to 0$ as $N\to +\infty$. Thus $\cesp{}{\Theta^N}\to 0$ as $N\to +\infty$. To conclude the first part of this proof, we can state that for any $q^0\in (0,+\infty)^V$, $\epsilon<\min_v q^0_v$, $T\leqslant \tau^\epsilon(q^*)$ and $\delta>0$,
		\[\cesp{}{\sup_{0 \leq t \leq T \wedge T_{\text{\textnormal{tube}}}^N}\norm{Q^N(t)-q^*(t)}{}}\to 0 \text{ as }N\to +\infty.\]
		
		Let us now remove the localization and prove that $Q^N \overset{\P}{\to} q^*$ uniformly on $[0,T]$. In order to do so, it is enough to show that $\P(T_{\text{\textnormal{tube}}}^N \geq T) \to 1$: indeed for any $\delta>0$,
		\begin{align*}
		\cpro{}{\sup_{t \leqslant T}\norm{Q^N(t)-q^*(t)}{\infty}\geqslant \delta}&\leqslant \cpro{}{\sup_{t \leqslant T\wedge T_{\text{\textnormal{tube}}}^N}\norm{Q^N(t)-q^*(t)}{\infty}\geqslant \delta, T_{\text{\textnormal{tube}}}^N\geqslant T }\\
		&\quad+\cpro{}{T_{\text{\textnormal{tube}}}^N< T}\\
		&\leqslant \cpro{}{\sup_{t \leqslant T\wedge T_{\text{\textnormal{tube}}}^N}\norm{Q^N(t)-q^*(t)}{\infty}\geqslant \delta}+\cpro{}{T_{\text{\textnormal{tube}}}^N< T}.
		\end{align*} 
		We already proved that the first term on the right hand side converges to 0, $\P(T_{\text{\textnormal{tube}}}^N \geq T) \to 1$ means that the second term also vanishes as $N\to +\infty$. By definition of $T_{\text{\textnormal{tube}}}^N$, we have
		
		\[ \left \lVert Q^N(T_{\text{\textnormal{tube}}}^N) - q^*(T_{\text{\textnormal{tube}}}^N) \right \rVert_1\ind_{T_{\text{\textnormal{tube}}}^N<+\infty} \geq \min(C_-,C_+/3)\ind_{T_{\text{\textnormal{tube}}}^N<+\infty}. \]
		
		Since $T_{\text{\textnormal{tube}}}^N \wedge T = T_{\text{\textnormal{tube}}}^N$ on the event $\{T_{\text{\textnormal{tube}}}^N \leq T\}$, this entails
		\[ \P \left( T_{\text{\textnormal{tube}}}^N \leq T \right) = \P \left( \left \lVert Q^N(T_{\text{\textnormal{tube}}}^N \wedge T) - q^*(T_{\text{\textnormal{tube}}}^N \wedge T) \right \rVert_1 \geq \min(C_-,C_+/3), \, T_{\text{\textnormal{tube}}}^N \leq T \right). \]
		Since we have proved that $Q^N(\, \cdot \, \wedge T_{\text{\textnormal{tube}}}^N)\overset{\P}{\to} q^*$ uniformly on $[0,T]$ for any $T\leqslant \tau^\epsilon(q^*)$, the previous probability vanishes.
		
		To sum up, for any $\epsilon,\delta>0$, for any $T\leqslant \tau^\epsilon(q^*)$, 
		\[\cpro{}{\sup_{0 \leqslant t \leqslant T}\norm{Q^N(t)-q^*(t)}{\infty}\geqslant \delta}\to 0 \text{ as } N\to +\infty.\]
		By continuity of $q^*$, since $\min_v q^0_v>0$, we have \[\lim_{\epsilon\to 0}\tau^\epsilon(q^*)= \tau^0(q^*)\] and \[q^*(\tau^\epsilon(q^*))\to q^*(\tau^0(q^*)) \text{ as } \epsilon\to 0.\] 
		For any $T<\tau^0(q^*)$, there is $\epsilon>0$ such that $T<\tau^\epsilon(q^*)$, and thus $Q^N\to q^*$ uniformly over compact time sets of $[0,\tau^0(q^*))$.

	\qed 
	\section{\large{Complete interference graph}}\label{sec:thm2}
	This section only deals with a complete interference graph so we omit dependency in the graph in all of the statements and note solutions to \eqref{eq:limitCIG} as $q^*$. Recall the definitions: 
	\[\tau^\epsilon(s\circ Q^N)=  \inf\lbrace t>0 ,\;  s(Q^N(t))\leqslant \epsilon\rbrace,\]
	and 
	\[\tau^\epsilon(s\circ q^*) =  \inf\lbrace t>0 ,\;  s(q^*)(t)\leqslant \epsilon\rbrace.\]
	Introduce as well
	\[T^{\epsilon}_-(f)\coloneqq \inf\lbrace t>0,\; \min_v f_v(t)\geqslant \epsilon\rbrace.\]
	We will prove that queue lengths become positive for small positive times and then identify the limit as the solution to \eqref{eq:limitCIG} using Theorem \ref{thm:main} up to $T_{\text{ext}}(q^0)$. This is where we have to distinguish three cases:
	
	\begin{itemize}
		\item (I) $s(\lambda)<1$,
		\item (II) $s(\lambda)> 1$,
		\item (III) $s(\lambda)=1$.
	\end{itemize}
	In (I) we will show that the sum of queue lengths stays absorbed at 0 when it reaches it. In case (II) we will show that queue lengths converge to $+\infty$ in all coordinates. For case (III) we will rely on the two previous cases: if $q^0=\textbf{0}$ the process stays absorbed at $0$, if $q^0\neq\textbf{0}$ all coordinates are positive in positive time and they do not reach $0$ in finite time. In the next section, we will prove Theorem \ref{thm:beyond} given a lemma governing the behaviour of queue lengths when some are near zero.
	\subsection{\large{\underline{Proof of Theorem \ref{thm:beyond}}}}
	
	Theorem \ref{thm:beyond} strengthens Theorem \ref{thm:main} to include any initial condition and stop at any finite time horizon. First we need a lemma roughly stating that it takes a time very small for the limit of $Q^N$ to cross small thresholds: 
	\begin{lemma}
		If $q^0\neq \textnormal{\textbf{0}}$ or $\sum_{v \in V}\lambda_v>1$, for any $\epsilon>0$ small enough, there is $t_\epsilon>0$ , such that \label{lem:pos}
		\begin{equation*}
		\cpro{q^0}{T_-^{\epsilon}(Q^N)\geqslant t_\epsilon}\to 0 \text{ as } N\to +\infty.
		\end{equation*}
		In addition, $t_\epsilon\to 0$ as $\epsilon\to 0$.
		
	\end{lemma}
	
\textbf{{\underline{Proof sketch:}}}
		
		The proof is based on an induction argument: we will recursively add nodes in a non-decreasing subset of nodes with bounded below queue lengths. To do that, we will define a succession of thresholds $\epsilon^v$, stopping times $\varphi^N_v$ and time horizons $t^v$ such that the set of nodes 
		\[J^N_i\coloneqq \lbrace v\in V ,\; Q^N_v(\varphi^N_i)>\epsilon^i\rbrace\]
		is increasing with $i$ and $\varphi^N_v\leqslant t^v$ almost surely. At step $k$ of the induction argument, we will prove that the probability of $\lbrace k\in J^N_k \rbrace$ goes to $1$ as $N\to +\infty$. The thresholds are defined in a way such that if $v$ is not in $J^N_{v-1}$, $Q^N_v$ increases on the mean, at least until the next threshold. On the other hand, once $v\in J^N_k$, we choose the time intervals such that for any $k'\geqslant k$, $v\in J^N_{k'}$ with high probability. Because the proof is quite technical and rely on a coupling argument, it has been separated in two parts. Section \ref{sec:mainstep}  will give the main steps of the proof and Section \ref{sec:coupling} will deal with the coupling argument. \qed 
	
	\begin{remark}
		Because of this result, we are able to extend Theorem 2.1 from \cite{us} for any initial condition $q^0\in \mathds{R}^V\setminus\lbrace\textbf{0}\rbrace$ with minimal work. Homogenization must hold at any positive time, and there is a strong attraction to the invariant manifold $I_\lambda$ so the result would be the same, with a possible jump at time zero to $I_\lambda$.
	\end{remark}
	Given this lemma, we just need to be careful when some initial states are null but the proof of Theorem \ref{thm:beyond} is straightforward using tightness and uniqueness of solutions of \eqref{eq:limitCIG}: the idea is that tightness is enough to ensure uniqueness of the limit by using uniqueness of solution to \eqref{eq:limitCIG}. The only other issue arises when some queue reach zero but it can only happen in a specific way that leaves all prelimit queues close to zero so it is consistent with the extension of maximal solution that we defined in \ref{def:ext}.\vspace{0.3cm}
	
\textbf{{\underline{Proof of Theorem \ref{thm:beyond}:}}}
		
		Let's sum up the proof method for this theorem.
		\begin{itemize}
			\item 	In case (I) the queue process reaches $\textbf{0}$ in finite time and stays absorbed.
			\item  In case (II), after some time all queues are positive and the dynamic evolves as in Theorem \ref{thm:main}. Starting from an initial condition positive at each node, Theorem \ref{thm:main} is enough to obtain a convergence uniform over compact sets.
			\item  	The case (III) can be handled similarly to case (I) or (II) depending on the initial state: the scaled sum remains constant as long as at least one queue is of order $N$. If $s(q^0)=0$, all components remain at $0$ as in case (I). Starting from a non trivial state, the queue lengths will evolve according to the ODE given by \eqref{eq:limitCIG}.  If $s(q^0)\neq 0$,  even if some queues start empty, all coordinates will be positive for positive time. Like for case (II), $\tau^0(q^*)=+\infty$ by Lemma \ref{lem:qpos}.\vspace{0.2cm}
		\end{itemize}
		
		By Proposition \ref{prop: ctight}, $Q^N$ is tight and any limiting point is a continuous function. Let $N_k$ be a subsequence such that $Q^{N_k}\Rightarrow \overline{q}$ as $k\to +\infty$ and $Q^{N_k}(0)\to q^0\in [0,+\infty)^V$. We know that $\overline{q}$ must be almost surely continuous by Proposition \ref{prop: ctight}.
		
		We first prove convergence of $Q^N$ up to $T_{\text{ext}}(q^0)$ starting from any initial condition. If $q^0_v>0$ for all $v\in V$, this is Theorem \ref{thm:main} so we assume $\min_v q^0_v=0$ for the first part of this proof. Since we already know that $Q^N$ is tight and any limiting point is in the space of continuous functions, we simply need to uniquely identify any potential limit.

		First, $\min_v\bar{q_v}$ is a continuous function almost surely because $\bar{q}$ is as well. Let us call $\bar{q}^{-1}$ the generalized inverse of $\min_v\bar{q_v}$: 
		\[\bar{q}^{-1}(t)\coloneqq \inf\lbrace s\geqslant 0,\; \min_v\bar{q_v}(s)\geqslant t\rbrace.\] 
		By Lemma 2.10 in Chapter 6 of \cite{jac03}, $\bar{q}^{-1}$ is almost surely càg and the set
		
		\[\mathcal{J}\coloneqq \lbrace t\geqslant 0,\; \cpro{}{\bar{q}^{-1}\text{ is not continuous in } t}>0\rbrace,\] of fixed time discontinuities of $\bar{q}^{-1}$ is countable. 
		Let $(\epsilon_p)_{p\in \entier}\in \mathcal{J}^c$ be a sequence of real number converging to 0. By definition of $\mathcal{J}$, $\bar{q}^{-1}$ is almost surely continuous at $\epsilon_p$ for every $p\in \entier$ and $\epsilon_p\to 0$ as $p\to +\infty$. Said otherwise, \[\cpro{}{\bar{q}^{-1} \text{continuous at } \epsilon_p}=1\] for all $p\in \entier$. Consider $T^\epsilon_-$ as an operator on càdlàg functions. By Proposition 2.11 of Chapter 6 in \cite{jac03}, $T^{\epsilon_p}_-$ is almost surely continuous at $\bar{q}$ because by definition of $\epsilon_p$, $\bar{q}^{-1}$ is continuous at $\epsilon_p$.
		
		By the continuous mapping theorem (Theorem 2.7 from \cite{bill}), since the limit of $Q^N$ is $\bar{q}$ and it is continuous and $T^{\epsilon_p}_-$ is almost surely continuous at $\bar{q}$, for any $p\in \entier$,
		\[(Q^{N_k}, T^{\epsilon_p}_-(Q^N))\Rightarrow (\bar{q}, T^{\epsilon_p}_-(\bar{q})).\]
		Using the joint convergence from the continuous mapping theorem, for any $p\in \entier$, 
		\[Q^{N_k}(T^{\epsilon_p}_-(Q^{N_k})+\cdot)\Rightarrow \bar{q}(T^{\epsilon_p}_-(\bar{q})+\cdot) \text{ as }k\to +\infty.\]
		Next, we explain how shifting the trajectory with the stopping time $T^{\epsilon_p}_-(Q^{N_k})$ gives another description of the limit using Theorem \ref{thm:main}.
		
		By definition of $T^{\epsilon_p}_-(Q^{N_k})$, for any $v\in V$, \[Q_v^{N_k}(T^{\epsilon_p}_-(Q^{N_k}))\geqslant \epsilon_p.\]
		We get $Q^{N_k}(T^{\epsilon_p}_-(Q^{N_k}))\Rightarrow \bar{q}(T^{\epsilon_p}_-(\bar{q}))$, with $\bar{q}_v(T^{\epsilon_p}_-(\bar{q}))\geqslant \epsilon_p$. By Theorem \ref{thm:main}, and the strong Markov property, as $k\to +\infty$,
		\[Q^{N_k}(T^{\epsilon_p}_-(Q^{N_k}+\cdot))\Rightarrow q^*(\cdot, \bar{q}(T^{\epsilon_p}_-(\bar{q}))).\]
		By uniqueness of the limit, it is necessary to have 
		\[ \bar{q}(T^{\epsilon_p}_-(\bar{q})+\cdot)=q^*(\cdot, \bar{q}(T^{\epsilon_p}_-(\bar{q}))).\]
		By Lemma \ref{lem:pos}, for any $p>0$, there is $t_p$ such that 
		\[\cpro{}{T_-^{\epsilon_p}(\bar{q})>t_p}=0,\]
		with $t_p\to 0$ as $p\to +\infty$.
		Since $T^{\epsilon_p}_-(\bar{q})\leqslant t_p$ almost surely and $t_p\to 0$ as $p\to +\infty$, $T^{\epsilon_p}_-(\bar{q})\to 0$ almost surely, and by continuity of $\bar{q}$, $\bar{q}(T^{\epsilon_p}_-(\bar{q}))\Rightarrow q^0$ as $p\to +\infty$. 
		
		In order to use uniform convergence to exchange the limits in $k$ and $p$, we use Skorohod representation theorem (see for instance Theorem 6.7 of \cite{bill}). It states that since $\reels_+^V$ is separable, and $Q^{N_k}\Rightarrow \bar{q}$ uniformly over compact sets, it is possible to construct a probability space such that $Q^{N_k}\to \bar{q}$ almost surely for the topology of uniform convergence on compact sets. Until the end of the proof, consider $Q^{N_k}$ and $\bar{q}$ constructed in such a way.

		To conclude, uniformly over compact sets of $[0,T_{\text{ext}}(q^0))$: for any $p$ and $N$ greater than $0$, 
		\begin{align*}
		\sup_{0 \leq t \leq T}\norm{Q^{N_k}(t)-q^*(t,q^0)}{\infty}&\leqslant \sup_{0 \leq t \leq T}\norm{Q^{N_k}(t)-Q^{N_k}(t+T^{\epsilon_p}(Q^N))}{\infty}\\
		&+\sup_{0 \leq t \leq T}\norm{Q^{N_k}(t+T_-^{\epsilon_p}(Q^{N_k}))-q^*(t,\bar{q}(T_-^{\epsilon_p}(\bar{q})))}{\infty}\\
		&+\sup_{0 \leq t \leq T}\norm{q^*(t,\bar{q}(T_-^{\epsilon_p}(\bar{q})))-q^*(t,q^0)}{\infty}.
		\end{align*}
		By letting $N\to +\infty$ then $p\to +\infty$, and using continuity moduli, we get
		\[Q^{N_k}\overset{\P}{\to} q^*(\cdot,q^0).\]
		
		In cases where $T_{\text{ext}}(q^0)=+\infty$, Theorem \ref{thm:beyond} is proved. 
		
		Together with the previous argument when queue start null, this regroups the case where $s(\lambda)\geqslant 1$. This also shows that in a complete interference graph with $s(\lambda)<1$, starting from any initial condition, $\tau^0(Q^N)\to \tau^0(s(q^*))$ in the mean. 
		
		We now focus on the case $s(\lambda)< 1$. Because of the expression of the sum in the limit, it will reach zero in finite time. We will now prove that when $Q^N(0)\to \textbf{0}$, the limit of the sum of queue lengths remains absorbed at $\textbf{0}$. Recall the definition of the stopping times $T^\epsilon_-(s\circ Q^N)$ given by
		
		\[T^{\epsilon}_-(s\circ Q^N)\coloneqq \inf\lbrace t>0 ,\; s(Q^N(t))\geqslant \epsilon\rbrace.\] 
		
		Here we emphasize the dependency in the trajectory but to ease notations later, we will only emphasize in the computation when it depends on the trajectory shifted after some time. The argument is essentially the same as with a deterministic system: to reach a level $\epsilon$, the sum should first reach $\epsilon/2$ and then be increasing between $\epsilon/2$ and $\epsilon$ which is highly unlickely.  Introduce the set $R^N(\epsilon)$ by 
		
		\[R^N(\epsilon)\coloneqq \lbrace q\in \frac{1}{N}\entier^V ,\; s(q)\in [\frac{\epsilon}{2}-\frac{1}{N}; \frac{\epsilon}{2}+\frac{1}{N})\rbrace.\]
		
		Note that for any $\epsilon>0$, $N<+\infty$,  $R^N(\epsilon)$ is a finite set. For $q\in \reels_+^V$, the notation $\P_q$ denotes $\P$ conditionally on $Q^N(0)=q$. Let's fix $\epsilon>0$ and call 
		
		\[f^N(t)\coloneqq \sup\limits_{q\in R^N(\epsilon)} \cpro{q}{T^{\epsilon}_-\leqslant t}.\]
		Since $R^N(\epsilon)$ is a finite set, there exists $q^N(t)\in R^N(\epsilon)$ such that $f^N(t)=\cpro{q^N(t)}{T^{\epsilon}_-\leqslant t}$.\\
		
		Fix $t<+\infty$ and $\epsilon>0$, assume $Q^N(0)\to \textbf{0}$ as $N\to +\infty$. We will prove that for $\epsilon>0$, the probability that the sum exceeds $\epsilon$ in finite time goes to $0$ as $N\to +\infty$. For any $\nu$ probability measure, we use $\P_\nu$ to denote the probability measure $\sum_{q\in \entier^V}\nu(q)\cpro{}{\cdot \mid Q(0)=q}$, if $\nu=\delta_q$ a Dirac measure, $\P_\nu=\P_q$. We get
		
		\begin{align}
		\cpro{Q^N(0)}{T^{\epsilon}_-\leqslant t}&=\notag\cpro{Q^N(0)}{T^{\epsilon/2}_-\leqslant t,T^{\epsilon}_-(s\circ(Q^N(\cdot+T^{\epsilon/2}_-))\leqslant t-T^{\epsilon/2}_-}\\
		&\notag\leqslant \cpro{Q^N(0)}{T^{\epsilon/2}_-\leqslant t,T^{\epsilon}_-(s\circ(Q^N(\cdot+T^{\epsilon/2}_-))\leqslant t}\\
		&= \cesp{Q^N(0)}{\ind_{T^{\epsilon/2}_-\leqslant t}\cpro{Q^N(T^{\epsilon/2}_-)}{T^{\epsilon}_-\leqslant t}},\label{eq:e/2beforee}
		\end{align}
		
		where the last equality is obtained using the strong Markov property. As it turns out, since jump size is $\frac{1}{N}$, almost surely on $\lbrace T^{\epsilon/2}_-\leqslant t\rbrace$, we have $Q^N(T^{\epsilon/2}_-)\in R^N(\epsilon)$ . For this reason, 
		
		\begin{equation}
		\cpro{Q^N(0)}{T^{\epsilon}_-\leqslant t}\leqslant \ind_{T^{c,\epsilon/2}_N\leqslant t}\cpro{Q^N(T^{c,\epsilon/2}_N)}{T^{\epsilon}_-\leqslant t}\leqslant f^N(t).
		\end{equation}
		
		We now prove that $f^N(t)\to 0$ for all $t$. Since $f^N$ is bounded by $1$, for any $t$ $f^N(t)$ is tight as a sequence of real numbers and we only need to prove uniqueness of the limit in terms of subsequence. For any $t$, there exists a subsequence such that $q^{N_k}(t)\to \widetilde{q}(t)$ because they too are tights. Any $\widetilde{q}(t)$ obtained in that way has $s(\widetilde{q}(t))=\frac{\epsilon}{2}$ because $q^N(t)\in R^N(\epsilon)$. \\
		By the first part of the proof, starting from any $q^N\in R^N(\epsilon)$ bounded away from $\textbf{0}$, the sum of queue lengths should decrease to $0$ before being able to reach $\epsilon$. Necessarily, for any initial condition in $ R^N(\epsilon)$, $\tau^0(s\circ Q^N)\overset{\P}{\to} \frac{\epsilon}{2(1-\sum_{v \in V}\lambda_v)}$ by the first part of the proof. Define $t^*\coloneqq \frac{\epsilon}{2(1-\sum_{v\in V}\lambda_v)}$. Assume that $t\geqslant t^*$, we get	
		\begin{align*}
		f^N(t)&=\cpro{q^N(t)}{T^{\epsilon}_-\leqslant t}\\
		&=\cpro{q^N(t)}{T^{\epsilon}_-\leqslant t, T^{\epsilon}_-\leqslant t^*}+\cpro{q^N(t)}{T^{\epsilon}_-\leqslant t, T^{\epsilon}_-> t^*}\\
		&\leqslant \cpro{q^N(t)}{T^{\epsilon}_-\leqslant t^*}+\cpro{q^N(t)}{T^{\epsilon}_-(s\circ Q^N(\cdot+t^*))\leqslant t-t^*}\\
		\end{align*}
		
		Similarly, by the Markov property,
		
		\begin{equation*}
		\cpro{q^N(t)}{T^{\epsilon}_-(Q^N(.+t^*))\leqslant t-t^*}= \cesp{q^N(t)}{\cpro{Q^N(t^*)}{T^{\epsilon}_-\leqslant t-t^*}}
		\end{equation*}
		
		Finally, since for any sequence $q^N\in R^N(\epsilon)$, under $\P_{q^N}$, $Q^N(t^*)\overset{\P}{\to} \textbf{0}$, we can apply the same reasoning as in \eqref{eq:e/2beforee} to prove that 	
		\[\cpro{Q^N(t^*)}{T^{\epsilon}_-\leqslant t-t^*}\leqslant f^N(t-t^*)\ind_{t\geqslant t^*}\]
		
		To sum this up, $f^N(t)\leqslant f^N(t^*)+f^N(t-t^*)\ind_{t\geqslant t^*}$, and iterating, we get 
		
		\[f^N(t)\leqslant \lceil \dfrac{t}{t^*}\rceil f^N( t^*).\]
		
		We can conclude the proof by using convergence up to $\tau^0(q^*)$ to prove that $f^N(t^*)\to 0$ as $N\to +\infty$: $f^N(t^*)=\cpro{q^N(t^*)}{T^\epsilon_-(s\circ Q^N)\leqslant t^*}$ but since the limiting sum is decreasing with $t$ and $q^N(t^*)\in R^N(\epsilon)$, we get $T_{\text{ext}}(q^N(t^*))\to t^*$ and the limiting process will reach \textbf{0} at $t^*$ because $s(q^N(t^*))\to \frac{\epsilon}{2}$.\\
	\qed

	\section{\large{Positivity of queue lengths for small times}}\label{app:c}
	
	The only thing left is to prove Lemma \ref{lem:pos}. It will be helpful to accurately describe a coupling for the different activation time/duration.  We present in the last section of the paper the proofs using an exact construction of the process.

	\subsection{\large{\underline{Explicit coupling}}}
	
	Let $\mathcal{E}(\lambda)$ denote the exponential distribution with parameter $\lambda$ and $\mathcal{G}(p)$ the geometric distribution with success probability $p$. We use the symbol $\overset{(d)}{=}$ to denote equality in distribution. Elementary computations show that if $G\overset{(d)}{=}\mathcal{G}(p)$ and $(E_k)_{k\in \entier}$ an \textit{i.i.d.} family of exponential variables with parameter $\lambda$ are independent,  \begin{equation*}
	\sum_{k=1}^{G}E_k\overset{(d)}{=}\mathcal{E}(p\lambda)\text{, }\min_{k=1,\ldots,n} E_k\overset{(d)}{=}\mathcal{E}(n\lambda)\text{ and } E_k\overset{(d)}{=}(\lambda)^{-1}\mathcal{E}(1).
	\end{equation*}
	Given the schedule, the dynamic of $Q$ is easy to construct. Arrivals at each node $v\in V$ happen at the points of independent Poisson processes on $\reels_+$ of intensities $\lambda_v$: let $(P_v)_{v\in V}$ be $n$ independent Poisson processes on $\reels_+$ of intensities $(\lambda_v)_{v\in V}$. They will serve as arrival processes. Similarly, let $(R_v)_{v\in V}$ also be independent unit intensity Poisson processes. Departures from queue $v\in V$ happen at the points of $\widetilde{R}_v(t)\coloneqq R_v(\int_{0}^{t}\ind_{Q_v(s)>0}\sigma_v(s)\texttt{d}s)$. All $P_v$ and $R_v$ are independent.
	
	For any initial condition $q^0$, 	introduce 
	\[V'\coloneqq\lbrace q\in \reels_+^V,\, q^0_v>0\rbrace\text{ and } \epsilon_0\coloneqq\min_{v\in V'}q^0_v.\]
	
	Lemma \ref{lem:pos} will be proved by induction on the nodes not in $V'$. From now on assume $V'\neq\emptyset$ and $V\setminus V'\neq \emptyset$. Fix $v^0\in V\setminus V'$ we will construct the dynamic between each activation of $v^0\in V\setminus V'$. 
	
	Let $(b(l))_{l\in \entier}$ and $(d(l))_{l\in \entier}$ the successive activation  and deactivation times of node $v^0$: assume $\sigma(0)=\textbf{0}$, then  $d(0)=0$ and for all $k>0$, 
	\begin{equation}
	\label{eq:tauannex}b(k)=\inf\{t>d(k-1) ,\, \sigma(t)=v^0\},\; d(k)=\inf\{t>b(k) ,\, \sigma(t)=\textbf{0}\} .
	\end{equation}

	Between each time the schedule is empty, $n$ inhomogeneous exponential variables compete to activate their nodes. For each $v\in V$, the activation rate is $\dfrac{1}{(Q_v(t)+1)^{-a}+1} $. \\
	Let's define $d^k_0=d(k)$, and for any $m\geqslant 1$,
	\[b^k_m\coloneqq \inf\left\lbrace t>d^{k}_{m-1},\, \sigma(t)\neq \textbf{0}\right\rbrace \text{ and } d^k_m\coloneqq \inf\left\lbrace t>b^k_{m},\, \sigma(s)=\textbf{0}\right\rbrace.\]
	More formally, let $(A^k_v(m))_{k\in \entier, m\in \entier, v \in V}$ be a family of unit parameter independent exponential variables. The $m^{th}$ time a node deactivates after $d(k)$, all nodes enter competition for the next activation. If this deactivation occurred before the $k+1^{th}$ activation of $v^0$, the deactivation happened at time $d^k_{m}$ and the next activation occurs at time \begin{equation}
	d^k_{m}+\min_{v\in V}\inf\left\lbrace t>0,\, A^{k}_v(m+1)\leqslant \int_{0}^{t}\dfrac{\d s}{1+(Q_v(d^k_{m}+s)+1)^{-a}}\right\rbrace,\label{activation}
	\end{equation}
	where for every $w$, before any activation,
	\begin{equation}
	Q_w(d^k_{m}+.)=Q_w(d^k_{m})+P_w(d^k_{m}+.)-P_w(d^k_{m}).\label{eq:evo}
	\end{equation}
	The queue that activates is the one realising the minimum in \eqref{activation}. Let $G(k)\geqslant 0$ be the number of times a node different from $v^0$ activates between the $k^{th}$ and $k+1^{th}$ activation of $v^0$. For any $k,m\in \entier$, the variables $(A^k_v(m))_{v\in V}$ are only used to determine the idle period and the active queue before the $m^{th}$ activation of a node after the $k^{th}$ activation of $v^0$ provided $v^0$ has not activated $k+1$ times. An important remark that we will use multiple times in this section is that the activation rates of any node can be bounded regardless of the state of the network: for any $Q\in \reels_+$,
	\[ 1/2\leqslant \dfrac{1}{1+(Q+1)^{-a}}\leqslant 1.\]
	
	Similarly, for any $m>0$, if $v$ is the $m^{th}$ node to activate after the $k^{th}$ activation of $v^0$, this activation occurs at time $b^k_m$, it will have an activation duration noted $d^k(m)-b^k(m)$ that can be expressed as another ``inhomogeneous exponential'' variable:  
	\begin{equation}
	d^k_m-b^k_m=\inf\left \lbrace t>0,\, D^k_m\leqslant \int_{0}^{t}\dfrac{\d s}{1+(Q_{v}(b^k_m+s)+1)^a}\right\rbrace,\label{deactivation}
	\end{equation}
	with $(D^k_m)_{k\in \entier, m\in \entier}$ unit intensity \textit{i.i.d.} exponential variables, discarded after each deactivation. For consistency purposes, the activation duration of node $v^0$ is constructed similarly with $(D^k_0)_{k\in \entier}$ a family of \textit{i.i.d.} exponential variables with unit parameter independent from the rest: for the $k^{th}$ activation period of $v^0$,
	\[d(k)-b(k)=\inf\left\lbrace t>0,\, D^k_0\leqslant \int_{0}^{t}\dfrac{\d s}{1+(Q_{v^0}(b(k)+s)+1)^a}\right\rbrace.\]
	
	Before any deactivation, the evolution of the active node is given by
	\[Q_v(b^k_m+t)=Q_v(b^k_m)+P_v(b^k_m+t)-P_v(b^k_m)-\widetilde{R}_v(b^k_m+t)+\widetilde{R}_v(b^k_m).\]
	Each inactive node evolves as in \eqref{eq:evo}.  Once node $v^0$ activates $k+1$ times, discard all $(A^k_v(m))_{v\in V, m\in \entier}$ and $(D^k_m)_{ m\in \entier}$. 
	
	\subsection{\large{\underline{General remarks}}}
	For all practical matters, it is possible to replace the Poisson processes for arrivals and departures by their ``compensated'' quantity by Lemma \ref{lem:mart}. Regardless of homogenization, for any $T<+\infty$ and $v\in V$,
	\begin{multline}
	\cpro{}{\sup_{t\leqslant T}\left\lvert Q_v^N(t)-Q_v^N(0)-\lambda_v t+\int_{0}^{t}\sigma_v^N(s)\ind_{Q_v^N(s)>0}\d s\right\rvert >N^{-1/3}}\\
	= \cpro{}{\sup_{t\leqslant T}\lvert M^N_v(t)\rvert> N^{-1/3}}\leqslant N^{1/3-1/2}\to 0 \text{ as } N\to +\infty.\label{eq:traj}
	\end{multline}
	
	In particular, since 
	\[\sup_{v\in V,\, t\leqslant T} Q_v^N(0)+\lambda_v t-\int_{0}^{t}\sigma^N_v(s)\ind_{Q^N_v(s)>0}\d s\leqslant \norm{Q^N(0)}{\infty}+\norm{\lambda}{\infty}T,\]
	we get that for any $v\in V$,
	\[\cpro{}{\sup_{t\leqslant T}Q_v^N(t)>\norm{q^0}{\infty}+\norm{\lambda}{\infty}T+N^{-1/3}}\to 0 \text{ as } N\to +\infty.\]
	
	All of the next discussion will be in the event 
	\[\left\lbrace \sup_{t\leqslant T} \lvert M^N(t)\rvert \leqslant N^{-1/3}\right\rbrace,\]
	and the $N^{-1/3}$ term will be omitted in the computations. This choice does not change the result and is only made to alleviate already cluttered notations. For the previous probability, we would for instance state that 
	\[\cpro{}{\sup_{t\leqslant T}Q_v^N(t)>\norm{q^0}{\infty}+\norm{\lambda}{\infty}T}\to 0 \text{ as } N\to +\infty.\]
	
	The proof of Lemma \ref{lem:pos} will be provided in a series of lemmas using a coupling argument. We will only deal with the case where $\norm{q^0}{\infty}>0$ because the case $s(\lambda)>1$ can be handled similarly: if $\norm{q^0}{\infty}=0$ and $s(\lambda)>1$, for any $t>0$ there is a queue bounded below. Since $s(\sigma^N(t))\leqslant 1$,  by summing over coordinates in \eqref{eq:traj},
	\[\cpro{}{s(Q^N(t))\geqslant s(Q^N(0))+(s(\lambda)-1)t-nN^{-1/3}}\to 1 \text{ as } N\to +\infty.\]
	In order for $s(q)$ to be greater than $\epsilon$, there is at least one $v$ such that $q_v\geqslant \frac{\epsilon}{n}$. Necessarily, with the abuse in notation explained above,
	\[\cpro{}{\max_{v\in V}Q^N_v(t)\geqslant \frac{(s(\lambda)-1)t}{n}}\to 1\text{ as }N\to +\infty.\]
	After any positive time, there is at least one queue with bounded below queue length. Recall the hitting time \vspace{0.2cm}
	\[T_-^{\epsilon}(Q^N_{v^0})=\inf\left\lbrace t>0,\, Q^N_{v^0}(t)\geqslant \epsilon\right\rbrace.\] In terms of hitting times, if $s(\lambda)>1$ and $\norm{q^0}{\infty}=0$,
	\[\cpro{}{T^\epsilon_-(\max_v Q^N_v)\leqslant \dfrac{n}{s(\lambda)-1}\epsilon}\to 1 \text{ as } N\to +\infty.\]
	With the strong Markov property, the reasoning in the next sections can be handled given which queue is above $\epsilon$ at time $T^\epsilon_-(\max_v Q^N_v)$ on the shifted process, and then integrate over $v\in V$.
	
	\subsection{\large{\underline{Main steps of the proof}}}\label{sec:mainstep}
Introduce
	\[K_{\textnormal{\text{time}}}\coloneqq \frac{8^{1+a}n\lVert q^0\rVert_{\infty}^a}{\min_v \lambda_v \epsilon_0^a}.\]
	
	Without loss in generality assume that $\min_v \lambda_v\leqslant 1$. Otherwise, on the fluid scale, all queue lengths are increasing and any level $\epsilon$ is trivially reached in a time linear in $\epsilon$. With the convention $1/0=+\infty$, let $\epsilon_1>0$ be such that
	\[\epsilon_1< \frac{\epsilon_0}{2}\min\left[1,\,\frac{1}{2^{1+1/a}}\left(\min_v \lambda_v \right)^{1/a},\,\dfrac{2\norm{q^0}{\infty}}{\norm{\lambda}{\infty}\epsilon_0K_{\textnormal{\text{time}}}},\,\dfrac{1}{K_{\textnormal{\text{time}}}(1-\min_v \lambda_v)}\right].\] 
	In this section, we first present the technical Lemma \ref{lem:iter} that allows us to prove Lemma \ref{lem:pos}. It provides a bound on crossing time of $\epsilon_1$ for $Q^N_{v^0}$. The probability that the hitting time exceeds a linear function of $\epsilon_1$ goes to $0$ as $N\to +\infty$. If the hitting time is small enough, queues that are in $V'$ do not have time decrease below $\epsilon_1$ before $Q^N_{v^0}$ increases above $\epsilon_1$.  After briefly justifying Lemma \ref{lem:iter}, we present the proof of Lemma \ref{lem:pos} using Lemma \ref{lem:iter}. The proof relies on an induction argument on the nodes in $V\setminus V'$. Next, we present some auxiliary lemmas for the proof of Lemma \ref{lem:iter}. Once the three auxiliary lemmas are stated the proof of Lemma \ref{lem:iter} is provided using a reasoning on events. Lemmas \ref{lem:notlong}, \ref{lem:proportion} and \ref{lem:nonneg} rely on a coupling argument.

	By definition of $\epsilon_0$,
	\[\min_{v\in V'}Q^N_v(t)\geqslant \epsilon_0-(1-\min_{v\in V}\lambda_v)t+M^N_{v}(t).\]
	Since
	\[\epsilon_1 K_{\textnormal{\text{time}}}<\dfrac{\epsilon_0}{2(1-\min_{v\in V}\lambda_v)},\]
	we get
	\begin{equation*}
	\cpro{}{\inf_{v\in V', s\leqslant \epsilon_1K_{\textnormal{\text{time}}}}Q^N_v(s)\geqslant \frac{\epsilon_0}{2}, \sup_{v\in V,\, s\leqslant \epsilon_1K_{\textnormal{\text{time}}}}Q^N_v(s)\leqslant \norm{q^0}{\infty}+\norm{\lambda}{\infty}\epsilon_1K_{\textnormal{\text{time}}}}\to 1 \text{ as } N\to +\infty.
	\end{equation*}
	Introduce the events
	\begin{equation}
	E_-^N(t)\coloneqq \left\lbrace \inf_{v\in V', s\leqslant t}Q^N_{v}(s)\geqslant  \frac{\epsilon_0}{2}\right\rbrace,\label{eq:event-}
	\end{equation}
	and 
	\begin{equation}
	E_+^N(t)\coloneqq \left\lbrace \sup_{s\leqslant t,\, v\in V}Q^N_v(s)\leqslant 2\norm{q^0}{\infty}\right\rbrace.\label{eq:event+}
	\end{equation}
	Since $\norm{\lambda}{\infty}\epsilon_1K_{\textnormal{\text{time}}}<\norm{q^0}{\infty}$,
	\begin{equation}
	\cpro{}{E_-^N(\epsilon_1K_{\textnormal{\text{time}}})\cap E_+^N(\epsilon_1K_{\textnormal{\text{time}}})}\to 1 \text{ as } N\to +\infty.\label{eq:bounddeter}
	\end{equation}

	\begin{lemma}\label{lem:iter}
		
		With the parameters defined before this lemma,
		\[\cpro{}{\min_{v\in V'\cup \lbrace v^0 \rbrace}Q_v^N(T_-^{\epsilon_1}(Q^N_{v^0}))\geqslant  \epsilon_1}\to 1 \text{ as } N\to +\infty,\]
		and 
		\[ \cpro{}{T_-^{\epsilon_1}(Q^N_{v^0})< \epsilon_1K_{\textnormal{\text{time}}}}\to 1 \text{ as } N\to +\infty.\]
	\end{lemma}
	Said with words, at the first time queue $v^0$ crossed $\epsilon_1$, all the other queues in $V'$ are also above $\epsilon_1$. The proof of this lemma relies entirely on the coupling presentedin the previous section. The general idea is to observe the evolution of node $v^0$ over $N^{1-a}K_{\textnormal{\text{act}}}$ activations of $v^0$. The proof of Lemma \ref{lem:iter} will be split in three auxiliary lemmas:
	\begin{itemize}
		\item The first lemma ensures that the large number of activations occurs with probability close to one in a time smaller than $N\epsilon_1K_{\textnormal{\text{time}}}$ so that a node that was in $V'$ does not have time to reach $0$ on the fluid scale. In addition bounding this time ensures that all queue lengths will remain bounded before it. This is necessary to be able to iterate the procedure and prove that the crossing time of a level $\epsilon$ converges to $0$ as $\epsilon\to 0$. 
		\item The second auxiliary lemma states that over this large number of activations of $v^0$, the fraction of the time this node was active is smaller than half its arrival rate, at least until it reaches $\epsilon_1$, provided all queues in $V'$ keep their queues lengths above $\frac{\epsilon_0}{2}$. If a queue is small compared to the other ones, it cannot be active for a large fraction of the time.
		\item  The third lemma states that this large number of activations occurs in a time of order of magnitude $N$. This means that over a non-negligible period of time queue $v^0$ will be increasing, at least until it reaches $\epsilon_1$.
	\end{itemize}
	Before formally stating these lemmas, we provide the proof of Lemma \ref{lem:pos} using Lemma \ref{lem:iter}.
	
\textbf{{\underline{Proof of Lemma \ref{lem:pos}:}}}

		By Lemma \ref{lem:iter}, 
		\[ \cpro{}{T_-^{\epsilon_1}(Q^N_{v^0})< \epsilon_1K_{\textnormal{\text{time}}}}\to 1 \text{ as } N\to +\infty,\]
		so $T^{\epsilon_1}_-(Q^N_{v^0})$ is a stopping time asymptotically finite. By the strong Markov property, the process $Q^N(\cdot+T_-^{\epsilon_1}(Q^N_{v^0}))$ is a Markov process with the same dynamic as $Q^N$ and starting point $Q^N(T_-^{\epsilon_1}(Q^N_{v^0}))$. Again by Lemma \ref{lem:iter}, we get
		\begin{equation}
		\cpro{}{\min_{v\in V'\cup \lbrace v^0 \rbrace}Q^N(T_-^{\epsilon_1}(Q^N_{v^0}))\geqslant  \epsilon_1}\to 1 \text{ as } N\to +\infty,\label{eq:iter1}
		\end{equation}
		and 
		\begin{equation}
		\label{eq:iter2}\cpro{}{\norm{Q^N(T^{\epsilon_1}_-(Q^N_{v^0}))}{\infty}>2\norm{q^0}{\infty}}\to 0 \text{ as } N\to +\infty,
		\end{equation}
		because $T^{\epsilon_1}_-(Q^N_{v^0}))\leqslant \epsilon_1K_{\textnormal{\text{time}}}$ and $\epsilon_1 K_{\textnormal{\text{time}}}<\frac{\norm{q^0}{\infty}}{\norm{\lambda}{\infty}}$.

		Recall the definition of $K_{\textnormal{\text{time}}}$:
		\[K_{\textnormal{\text{time}}}\coloneqq \frac{8^{1+a}n\lVert q^0\rVert_{\infty}^a}{\min_v \lambda_v \epsilon_0^a}.\]
		
		As long as $V\setminus V'\cup \lbrace v^0\rbrace \neq \emptyset$, because of \eqref{eq:iter1} and \eqref{eq:iter2}, it is possible to iterate the procedure with $V'$ replaced by $V'\cup \lbrace v^0\rbrace$, $\norm{q^0}{\infty}$ replaced by  $2\norm{q^0}{\infty}$ and $\epsilon_0$ replaced by $\epsilon_1$ for the shifted process $Q^N(\cdot+T_-^{\epsilon_1}(Q^N_{v^0}))$.
		
		When iterating the procedure at step $k>1$,  introduce the constant:
		\[K_{\textnormal{\text{time}}}(k)\coloneqq \frac{8^{1+a}n2^{a(k-1)}\lVert q^0\rVert_{\infty}^a}{\min_v \lambda_v \epsilon_{k-1}^a}.\]
		
		\[\epsilon_{k+1}< \frac{\epsilon_k}{2}\min\left[1,\,\frac{1}{2^{1+1/a}}\left(\min_v \lambda_v \right)^{1/a},\,\dfrac{2^{k}\norm{q^0}{\infty}}{\norm{\lambda}{\infty}\epsilon_kK_{\textnormal{\text{time}}}(k)},\,\dfrac{1}{K_{\textnormal{\text{time}}}(k)(1-\min_v \lambda_v)}\right],\] 
		and $\epsilon\leqslant\epsilon_{\lvert V\setminus V'\rvert} $, which is possible for $\epsilon$ small enough.
		In this case,
		\[\cpro{}{T^\epsilon_-(Q^N)\geqslant \left(\sum_{k=1}^{\lvert V\setminus V'\rvert}K_{\textnormal{\text{time}}}(k)\epsilon_k\right)}\to 0 \text{ as } N\to +\infty.\]
		When $\epsilon\to 0$, it is possible to choose $\sum_{k=1}^{\lvert V\setminus V'\rvert}\epsilon_k^{1-a}$ as small as desired.
	\qed  
	We now present the three lemmas used to prove Lemma \ref{lem:iter}. Their proofs are provided in the next section. Once those three lemmas are formally stated, we present a proof of Lemma \ref{lem:iter} at the end of the section. Introduce the constant:
	\[K_{\textnormal{\text{act}}}\coloneqq \dfrac{\epsilon_1K_{\textnormal{\text{time}}}}{n4^{1+a}\norm{q^0}{\infty}^a},\]
	and the stopping time
	\[d(N^{1-a}K_{\textnormal{\text{act}}})=\text{ Time of the }N^{1-a}K_{\textnormal{\text{act}}}^{th} \text{ deactivation for } v^0,\]
	see \eqref{eq:tauannex} for more details.
	\begin{lemma}\label{lem:notlong}

		As  $N\to +\infty$,
		\[\cpro{}{d(N^{1-a}K_{\textnormal{\text{act}}}))\leqslant N\epsilon_1K_{\textnormal{\text{time}}}}\to 1.\]
		Consequently,
		\[\cpro{}{E^N_-(\frac{d(N^{1-a}K_{\textnormal{\text{act}}})}{N})}=\cpro{}{\inf_{v\in V', Ns\leqslant d(N^{1-a}K_{\textnormal{\text{act}}}))}Q^N_{v}(s)\geqslant  \frac{\epsilon_0}{2}}\to 1 \text{ as } N \to +\infty,\]
		and 
		\[\cpro{}{E^N_+(\frac{d(N^{1-a}K_{\textnormal{\text{act}}})}{N})}=\cpro{}{ \sup_{v\in V,\, Ns\leqslant d(N^{1-a}K_{\textnormal{\text{act}}}))}Q^N_v(s)\leqslant 2\norm{q^0}{\infty}}\to 1 \text{ as } N \to +\infty.\]\vspace{0.3cm}
	\end{lemma}
	The idea is to reason in the event 
	\[E^N_+(\epsilon_1K_{\textnormal{\text{time}}})=\left\lbrace  \sup_{v\in V,\, s\leqslant \epsilon_1K_{\textnormal{\text{time}}}}Q^N_v(s)\leqslant 2\norm{q^0}{\infty}\right\rbrace,\]
	and use the fact that its probability goes to 1.
	When proving that $d(N^{1-a}K_{\textnormal{\text{act}}})$ is smaller than $N\epsilon_1 K_{\textnormal{\text{time}}}$, only the value of queue lengths before $N\epsilon_1 K_{\textnormal{\text{time}}}$ is important. In $E^N_+(\epsilon_1K_{\textnormal{\text{time}}})$, it is possible to bound each activation duration using the coupling described in the next section. The second part of the result is a direct consequence of \eqref{eq:bounddeter}.

	Introduce the event\begin{equation}\label{eq:eventap}
	E^N_{v^0}\coloneqq \left\lbrace T^{\epsilon_1}_-(Q^N_{v^0})> d(N^{1-a}K_{\textnormal{\text{act}}})\right\rbrace\cap E^N_-(d(N^{1-a}K_{\textnormal{\text{act}}})).
	\end{equation}
	The next lemma uses this event:
	\begin{lemma}\label{lem:proportion}
		As $N\to +\infty$, we get
		\[\cpro{}{\int_{0}^{d(N^{1-a}K_{\textnormal{\text{act}}})}\sigma_{v^0}(s)\d s> d(N^{1-a}K_{\textnormal{\text{act}}})\dfrac{\lambda_{v^0}}{2}, E^N_{v^0}}\to 0.\]
		In particular, because of Lemma \ref{lem:notlong} 
		\[\cpro{}{T^{\epsilon_1}_-(Q^N_{v^0})\leqslant d(N^{1-a}K_{\textnormal{\text{act}}})\text{ or } Q^N_{v^0}(d(N^{1-a}K_{\textnormal{\text{act}}}))\geqslant \frac{\lambda_{v^0}}{2}d(N^{1-a}K_{\textnormal{\text{act}}})}\to 1 \text{ as } N\to +\infty.\]
	\end{lemma}
	This is quite intuitive: as long as queue $v^0$ is small enough and all queues in $V'$ are large, queue $v^0$ cannot be active for a large fraction of the time compared to nodes in $V'$. This entails that $Q^N_{v^0}$ must be increasing during this period until it reaches $\epsilon_1$ with probability close to one.

	\begin{lemma}\label{lem:nonneg}
		With the constants described at the beginning of the section,
		\[\cpro{}{d(N^{1-a}K_{\textnormal{\text{act}}})<\frac{2N\epsilon_1}{\lambda_{v^0}},\, E^N_-(d(N^{1-a}K_{\textnormal{\text{act}}}))}\to 0 \text{ as } N\to +\infty.\]
		In particular, because of Lemma \ref{lem:notlong}, 
		\[\cpro{}{d(N^{1-a}K_{\textnormal{\text{act}}})\geqslant \frac{2N\epsilon_1}{\lambda_{v^0}}}\to 1 \text{ as } N\to +\infty.\]
	\end{lemma} 
	
	Given these three lemmas, the proof of Lemma \ref{lem:iter} can be based on the study of events whose probability goes to one:\vspace{0.3cm}
	
\textbf{{\underline{Proof of Lemma \ref{lem:iter}:}}}

		In terms of event,
		\begin{equation}\label{eq:event1}
		\left\lbrace T^{\epsilon_1}_-(Q^N_{v^0})\leqslant d(N^{1-a}K_{\textnormal{\text{act}}})\right\rbrace\cap \left\lbrace \frac{2N\epsilon_1}{\lambda_{v^0}}\leqslant d(N^{1-a}K_{\textnormal{\text{act}}})\leqslant N\epsilon_1K_{\textnormal{\text{time}}}\right\rbrace\subset \left\lbrace T^{\epsilon_1}_-(Q^N_{v^0})\leqslant N\epsilon_1K_{\textnormal{\text{time}}}\right\rbrace.
		\end{equation}
		Similarly,
		
		\[\left\lbrace\int_{0}^{d(N^{1-a}K_{\textnormal{\text{act}}})}\sigma_{v^0}(s)\d s\leqslant  d(N^{1-a}K_{\textnormal{\text{act}}})\dfrac{\lambda_{v^0}}{2}\right\rbrace\subset \left\lbrace Q^N_{v^0}(\frac{d(N^{1-a}K_{\textnormal{\text{act}}})}{N})\geqslant \dfrac{\lambda_{v^0}d(N^{1-a}K_{\textnormal{\text{act}}})}{2N}\right\rbrace,\]
		and
		
		\begin{multline}
		\left\lbrace Q^N_{v^0}(\frac{d(N^{1-a}K_{\textnormal{\text{act}}})}{N})\geqslant \dfrac{\lambda_{v^0}d(N^{1-a}K_{\textnormal{\text{act}}})}{2N}\right\rbrace\cap \left\lbrace \frac{2N\epsilon_1}{\lambda_{v^0}}\leqslant d(N^{1-a}K_{\textnormal{\text{act}}})\leqslant N\epsilon_1K_{\textnormal{\text{time}}}\right\rbrace\\
		\subset\left\lbrace T^{\epsilon_1}_-(Q^N_{v^0})\leqslant N\epsilon_1K_{\textnormal{\text{time}}}\right\rbrace\label{eq:event2}.
		\end{multline}
		
		By \eqref{eq:event1} and \eqref{eq:event2}, $\left\lbrace T^{\epsilon_1}_-(Q^N_{v^0})\leqslant N\epsilon_1K_{\textnormal{\text{time}}}\right\rbrace$ contains 
		\begin{multline*}\left(\left\lbrace NT^{\epsilon_1(Q^N_{v^0})}\leqslant d(N^{1-a}K_{\textnormal{\text{act}}})\right\rbrace
		\cup \left\lbrace Q^N_{v^0}(\frac{d(N^{1-a}K_{\textnormal{\text{act}}})}{N})\geqslant \dfrac{\lambda_{v^0}d(N^{1-a}K_{\textnormal{\text{act}}})}{2N}\right\rbrace\right)\\
		\cap\left\lbrace \frac{2N\epsilon_1}{\lambda_{v^0}}\leqslant d(N^{1-a}K_{\textnormal{\text{act}}})\leqslant N\epsilon_1K_{\textnormal{\text{time}}}\right\rbrace. \end{multline*}
		In conclusion,
		\begin{multline*}
		\cpro{}{T^{\epsilon_1}(Q^N_{v^0})\leqslant N\epsilon_1K_{\textnormal{\text{time}}}}\geqslant \mathds{P}\big(\left\lbrace \frac{2N\epsilon_1}{\lambda_{v^0}}\leqslant d(N^{1-a}K_{\textnormal{\text{act}}})\leqslant N\epsilon_1K_{\textnormal{\text{time}}}\right\rbrace\\
		\cap\left(\left\lbrace NT^{\epsilon_1(Q^N_{v^0})}\leqslant d(N^{1-a}K_{\textnormal{\text{act}}})\right\rbrace
		\cup \left\lbrace\int_{0}^{d(N^{1-a}K_{\textnormal{\text{act}}})}\sigma_{v^0}(s)\d s\leqslant  d(N^{1-a}K_{\textnormal{\text{act}}})\dfrac{\lambda_{v^0}}{2}\right\rbrace\right)\big).
		\end{multline*}
		We now justify why the probability on the right hand side converges to 1.	As we already mentioned, because of Lemma \ref{lem:proportion},
		\begin{equation*}
		\cpro{}{\left\lbrace NT^{\epsilon_1}_-(Q^N_{v^0})\leqslant d(N^{1-a}K_{\textnormal{\text{act}}})\right\rbrace\cup \left\lbrace Q^N_{v^0}(\frac{d(N^{1-a}K_{\textnormal{\text{act}}})}{N})\geqslant \dfrac{\lambda_{v^0}d(N^{1-a}K_{\textnormal{\text{act}}})}{2N}\right\rbrace}\to 1 \text{ as } N\to +\infty.
		\end{equation*}
		By Lemmas \ref{lem:notlong} and \ref{lem:nonneg}, 
		\[\cpro{}{\frac{2N\epsilon_1}{\lambda_{v^0}}\leqslant d(N^{1-a}K_{\textnormal{\text{act}}})\leqslant N\epsilon_1K_{\textnormal{\text{time}}}}\to 1 \text{ as } N\to +\infty.\]
		Since both probabilities converge to 1, the probability of the intersection converges to 1 as well.
		
		For the second part, in the event 
		\[\left\lbrace T^{\epsilon_1}_-(Q^N_{v^0})\leqslant \epsilon_1K_{\textnormal{\text{time}}}\right\rbrace,\]
		we get 
		\[\inf_{v\in V', s\leqslant T^{\epsilon_1}_-(Q^N_{v^0})}Q^N_v(s)\geqslant \inf_{v\in V', s\leqslant \epsilon_1K_{\textnormal{\text{time}}}}Q^N_v(s).\]
		Since $\epsilon_1<\frac{\epsilon_0}{2}$ and $\cpro{}{E^N_-(\epsilon_1K_{\textnormal{\text{time}}})}\to 1$, we get by definition of $E^N_-(t)$ that
		\[\cpro{}{\inf_{v\in V', s\leqslant  T^{\epsilon_1}_-(Q^N_{v^0})}Q^N_v(s)\geqslant \epsilon_1}\to 1 \text{ as } N\to +\infty.\]
		By definition of the stopping time, $Q^N_{v^0}(T^{\epsilon_1}_-(Q^N_{v^0}))\geqslant \epsilon_1$, which proves the second part of the result. 
	\qed 
	
	\subsection{\large{\underline{Proofs using coupling}}}\label{sec:coupling}
	
\textbf{{\underline{Proof of Lemma \ref{lem:notlong}}}}

		The first step is to consider the intersection between $E^N_+(\epsilon_1 K_{\textnormal{\text{time}}})$ and $\left\lbrace d(N^{1-a}K_{\textnormal{\text{act}}})\leqslant N\epsilon_1 K_{\textnormal{\text{time}}}\right\rbrace $.
		
		In $E^N_+(\epsilon_1 K_{\textnormal{\text{time}}})$, if the $m^{th}$ activation of a node between the $k^{th}$ and $k+1^{th}$ activations  of node $v^0$ is such that 
		\[d^k_m\leqslant N\epsilon_1K_{\textnormal{\text{time}}},\]
		by \eqref{deactivation}, $d_m^k-b_m^k$ is smaller than 
		\[\inf\left\lbrace t>0,\, D^k_m\leqslant \dfrac{t}{(N4\norm{q^0}{\infty})^{a}}\right\rbrace=(N4\norm{q^0}{\infty})^{a}D^k_m.\] 
		Similarly, for any $k\leqslant N^{1-a}K_{\textnormal{\text{act}}}$,
		\[d(k)-b(k)\leqslant \inf\left\lbrace t>0,\, D^k_0\leqslant \dfrac{t}{(N4\norm{q^0}{\infty})^{a}}\right\rbrace=(N4\norm{q^0}{\infty})^{a}D^k_0.\] 
		
		There are exactly $N^{1-a}K_{\textnormal{\text{act}}}$ activations of node $v^0$ before $d(N^{1-a}K_{\textnormal{\text{act}}})$. Between each activation of node $v^0$, the number of times a node other than $v^0$ activates can be bounded by a geometric variable with parameter $\frac{1}{2n-1}$ because of \eqref{activation}: let's call $G(k)$ the number of times a node different from $v^0$ activates between the $k^{th}$ and $k+1^{th}$ activation of $v^0$. By \eqref{activation}, $G(k)$ is smaller than the number of times a queue different from $v^0$ would activate if the activation rate of node $v^0$ was $\frac{1}{2}$ and the activation rates of all other queues were $1$. More formally, $G(k)$ is smaller than $\bar{G}(k)$: the number of $p>0$ such that
		\[\min\left[\min_{v\neq v^0}A^k_v(p), 2A^k_{v^0}(p)\right]\neq 2A^k_{v^0}(p),\]
		before the first time
		\[\min\left[\min_{v\neq v^0}A^k_v(p), 2A^k_{v^0}(p)\right]= 2A^k_{v^0}(p).\]
		Queue $v^0$ will necessarily activate if  \[\min\left[\min_{v\neq v^0}A^k_v(p), 2A^k_{v^0}(p)\right]= 2A^k_{v^0}(p),\]
		but the activation of $v^0$ may occur for a smaller $p'$, thus $G(k)\leqslant \bar{G}(k)$.
		
		By independence of $(A^m_v(p))_{m,p\in \entier, v\in V}$, \begin{equation}
		\bar{G}(k)\overset{(d)}{=}\mathcal{G}(\frac{1}{2n-1})\label{eq:geom}
		\end{equation} 
		
		and is independent from $(D^k_v)_{k\in \entier, v\in V}$. For any $p>0$, the queue realizing the minimum of the $p^{th}$ competition is independent from $\min\left[\min_{v\neq v^0}A^k_v(p), 2A^k_{v^0}(p)\right]$ and from the result of the competition for a different $k$ or $p$. This is due to independence of $(A^k_v(p))_{k,p\in \entier, v\in V}$ and the fact that for any $v\neq v^0$,
		\begin{equation}
		\cpro{}{\min\left[\min_{v\neq v^0}A^k_v(p), 2A^k_{v^0}(p)\right]=A^k_{v^1},\, A^k_{v^1}(p)>t}=\dfrac{2}{2n-1}\exp(-t\dfrac{2n-1}{2}) ,\label{eq:geom1}
		\end{equation}
		and 
		\begin{equation}
		\cpro{}{\min\left[\min_{v\neq v^0}A^k_v(p), 2A^k_{v^0}(p)\right]=2A^k_{v^0}(p),\, 2A^k_{v^0}(p)>t}=\frac{1}{2n-1}\exp(-(\frac{2n-1}{2})t).\label{eq:geom2}
		\end{equation}
		
		Since the activation rate is greater than $\frac{1}{2}$, the associated idle period is smaller than \[\widetilde{A}^k(p)\coloneqq 2\min\left[\min_{v\neq v^0}A^k_v(p), 2A^k_{v^0}(p)\right]\overset{(d)}{=}\frac{4}{2n-1}\mathcal{E}(1),\] which is independent from $\bar{G}(k)$. The independence between $\bar{G}(k)$ and the duration of each idle period comes from \eqref{eq:geom1}, \eqref{eq:geom2} and independence of $(A^k_v(p))_{k,p\in \entier, v\in V}$. Let us define

		\[ \bar{D}^N\coloneqq \sum_{k=1}^{N^{1-a}K_{\textnormal{\text{act}}}}\sum_{p=0}^{\bar{G}(k)}\left(\widetilde{A}^k(p)+(4N\norm{q^0}{\infty})^a D^k_p\right),\]
		with the convention that $\widetilde{A}^k_0\coloneqq \widetilde{A}^k_{\bar{G}(k)}$, which is still independent from $(\bar{G}(k))_{k\in \entier}$. The term with $p=0$ handles the activation duration of node $v^0$ and the last idle period before $v^0$ activates. Because of the use of different $(A^k_v(m))_{v\in V, m\in \entier}$ after the $k+1^{th}$ activation of $v^0$, $(\bar{G}(k))_{k\in \entier}$ forms an \textit{i.i.d.} family of geometric variables. Hence, $\bar{D}^N$ is a sum of \textit{i.i.d.} variables.

		Recall the event \[E^N_+(\epsilon_1K_{\textnormal{\text{time}}})=\left\lbrace \sup_{v\in V,\, s\leqslant \epsilon_1K_{\textnormal{\text{time}}}}Q^N_v(s)\leqslant 2\norm{q^0}{\infty}\right\rbrace.\] 
		The coupling ensures that 
		\begin{align*}
		\cpro{}{d(N^{1-a}K_{\textnormal{\text{act}}})\leqslant N\epsilon_1K_{\textnormal{\text{time}}},\, E^N_+(\epsilon_1K_{\textnormal{\text{time}}})}&\geqslant \cpro{}{ \bar{D}^N\leqslant N\epsilon_1K_{\textnormal{\text{time}}},\, E^N_+(\epsilon_1K_{\textnormal{\text{time}}})}\\
		&=\cpro{}{\frac{1}{NK_{\textnormal{\text{act}}}} \bar{D}^N\leqslant \dfrac{\epsilon_1K_{\textnormal{\text{time}}}}{K_{\textnormal{\text{act}}}},\, E^N_+(\epsilon_1K_{\textnormal{\text{time}}})}
		\end{align*}

		By construction of the coupling, $(\widetilde{A}^k(p))_{k,p\in \entier}$ and $(D^k_p)_{k,p\in \entier}$ from two \textit{i.i.d.} families, independent from $\bar{G}(k)$ for any $k>0$. By independence of $(\bar{G}(k))_{k\in \entier}$, $(\widetilde{A}^k_p)_{k,p\in \entier}$,$ (D^k_p)_{k,p\in\entier}$ and independence between those variables, by the law of large numbers 
		
		\begin{multline}\dfrac{1}{N^{1-a}K_{\textnormal{\text{act}}}}\sum_{k=1}^{N^{1-a}K_{\textnormal{\text{act}}}}\sum_{p=0}^{\bar{G}(k)}\left(\frac{\widetilde{A}^k_p}{N^a}+(4\norm{q^0}{\infty})^a D^k_p\right)\to n2^{1+2a}\norm{q^0}{\infty}^a \text{ almost surely as } N\to +\infty.\end{multline}
		By definition of $K_{\textnormal{\text{act}}}$,	\[K_{\textnormal{\text{act}}}=\dfrac{\epsilon_1K_{\textnormal{\text{time}}}}{4^{1+a}n\norm{q^0}{\infty}^a}<\dfrac{\epsilon_1K_{\textnormal{\text{time}}}}{2^{1+2a}n\norm{q^0}{\infty}^a},\]
		and so
		\[n2^{1+2a}\norm{q^0}{\infty}^a<\dfrac{\epsilon_1 K_{\textnormal{\text{time}}}}{K_{\textnormal{\text{act}}}}.\]

		From \eqref{eq:bounddeter}, we get 
		\[\cpro{}{E^N_+(N\epsilon_1K_{\textnormal{\text{time}}})}\to 1 \text{ as } N\to +\infty,\]
		
		and so,
		\[\cpro{}{d(N^{1-a}K_{\textnormal{\text{act}}})\leqslant N\epsilon_1K_{\textnormal{\text{time}}}}\to 1 \text{ as } N\to +\infty.\]
		The second part of the result comes from the fact that on $\lbrace d(N^{1-a}K_{\textnormal{\text{act}}})\leqslant N\epsilon_1K_{\textnormal{\text{time}}}\rbrace$, we get
		\[\sup_{v\in V,\,t\leqslant d(N^{1-a}K_{\textnormal{\text{act}}})}Q^N_v(t)\leqslant \sup_{v\in V,\, t\leqslant N\epsilon_1K_{\textnormal{\text{act}}}}Q^N_v(t)\leqslant 2\norm{q^0}{\infty},\]
		and 
		\[\inf_{v\in V',\,t\leqslant d(N^{1-a}K_{\textnormal{\text{act}}})}Q^N_v(t)\geqslant \inf_{v\in V',\, t\leqslant N\epsilon_1K_{\textnormal{\text{act}}}}Q^N_v(t)\geqslant \frac{\epsilon_0}{2}.\]
		In this situation, \eqref{eq:bounddeter} and the definition of the events in \eqref{eq:event-} and \eqref{eq:event+} are sufficient to state that 
		\[\cpro{}{E^N_+(\frac{d(N^{1-a}K_{\textnormal{\text{act}}})}{N})\cap E^N_-(\frac{d(N^{1-a}K_{\textnormal{\text{act}}})}{N})\cap \lbrace d(N^{1-a}K_{\textnormal{\text{act}}})\leqslant N\epsilon_1K_{\textnormal{\text{time}}}\rbrace}\to 1 \text{ as } N\to +\infty.\]
	\qed 
	We now turn to the proof of Lemma \ref{lem:proportion}. Recall the definition of the event: 
	\[E^N_{v^0}=\left\lbrace \sup_{Ns\leqslant d(N^{1-a}K_{\textnormal{\text{act}}})}Q^N_{v^0}(s)\leqslant \epsilon_1\,\right\rbrace\cap E^N_-( d(N^{1-a}K_{\textnormal{\text{act}}})).\]\vspace{0.3cm}
	
\textbf{{\underline{Proof of Lemma \ref{lem:proportion}:}}}

		By definition of the activation periods, 
		\[\int_{0}^{d(N^{1-a}K_{\textnormal{\text{act}}})}\sigma_{v^0}(s)\d s=\sum_{k=1}^{N^{1-a}K_{\textnormal{\text{act}}}}\left(d(k)-b(k)\right).\]
		
		In  $E^N_{v^0}$, for any $k\leqslant N^{1-a}K_{\textnormal{\text{act}}}$, 
		\[d(k)-b(k)\leqslant \inf\left\lbrace t>0,\, D^k_0\leqslant \frac{t}{(2\epsilon_1)^a}\right\rbrace=(2\epsilon_1)^a D^k_0,\]
		by \eqref{deactivation}, because in the event $E^N_{v^0}$, $Q^N_{v^0}(t)$ is bounded above by $\epsilon_1$ for $t\leqslant d(N^{1-a}K_{\textnormal{\text{act}}})$. By independence of $(D^k_v)_{k\in \entier, v\in V}$ and the deterministic nature of the parameter, the bound on activation durations form an \textit{i.i.d.} family of exponential variables. Similarly, it is possible to lower bound $d(N^{1-a}K)$ using once again \eqref{deactivation}. In $E^N_{v^0}$, every node in $V'$ has queue lengths uniformly bounded below by $\frac{\epsilon_0}{2}$ until $d(N^{1-a}K_{\textnormal{\text{act}}})$. For any $k<N^{1-a}K_{\textnormal{\text{act}}}$ and $m<G(k)$, such that $\sigma_v(b^k_m)=1$ for some $v\in V'$,
		\begin{equation*}
		d^k_m-b^k_m\geqslant (\frac{N\epsilon_0}{2})^aD^k_m.
		\end{equation*}
		
		There are $N^{1-a}K_{\textnormal{\text{act}}}$ activations for node $v^0$ before $d(N^{1-a}K_{\textnormal{\text{act}}})$. Inbetween activations for node $v^0$, there is a number of times where a node in $V'$ activates. Because of \eqref{activation}, the activation rate of a node in $V'$ is greater than $\frac{1}{2}$. The maximum activation rate of any queue is $1$. After any activation, the probability that a node in $V'$ is the next to activate is greater than $\frac{\lvert V'\rvert}{2n-\lvert V'\rvert}$ regardless of the evolution of the network, each try being independent from the others because of the competition between new $A^k_v(m)$.  The number of times a queue in $V'$ activates between the $k^{th}$ and $k+1^{th}$ activations of $v^0$ is greater than $\check{G}(k)$ constructed similarly to \eqref{eq:geom}: a geometric variable counting the number of times a queue in $V'$ activates before a queue not in $V'$ activates:
		\[\check{G}(k)\overset{(d)}{=}\mathcal{G}(\frac{\lvert V'\rvert}{2n-\lvert V'\rvert}).\]
		The independence between $(A^k_v(m))_{k,m\in \entier, \, v\in V}$ and $(D^k_m)_{k,m\in \entier}$ and deterministic nature of the intensities ensures that $\check{G}(k)$ is independent from $(D^k_m)_{k,m\in \entier}$. This discussion amounts to the bounds 
		\begin{equation*}
		\int_{0}^{d(N^{1-a}K_{\textnormal{\text{act}}})}\sigma_{v^0}(s)\d s\leqslant(2N\epsilon_1)^a\sum_{k=1}^{N^{1-a}K_{\textnormal{\text{act}}}}D^k_0, 
		\end{equation*}
		and
		\begin{equation}\label{eq:lowerd}
		d(N^{1-a}K_{\textnormal{\text{act}}})\geqslant \sum_{k=1}^{N^{1-a}K_{\textnormal{\text{act}}}}\sum_{m=1}^{\check{G}(k)}(\frac{N\epsilon_0}{2})^a D^k_m.
		\end{equation}
		Let's call 
		\[\overline{D}_k=\sum_{m=1}^{\check{G}(k)}(\frac{N\epsilon_0}{2})^a D^k_m.\]
		In the event $E^N_{v^0}$, we get
		\begin{align*}
		\cpro{}{\sum_{k=1}^{N^{1-a}K_{\textnormal{\text{act}}}}d(k)-d(k)> d_v(N^{1-a}K_{\textnormal{\text{act}}})\dfrac{\lambda_{v^0}}{2},\, E^N_{v^0}}&\leqslant  \cpro{}{(2N\epsilon_1)^a\sum_{k=1}^{N^{1-a}K_{\textnormal{\text{act}}}}D^k_0> \sum_{k=1}^{N^{1-a}K_{\textnormal{\text{act}}}}\overline{D}_k\dfrac{\lambda_{v^0}}{2},\,E^N_{v^0}}\\
		&\leqslant \cpro{}{\dfrac{1}{N^{1-a}K_{\textnormal{\text{act}}}}\sum_{k=1}^{N^{1-a}K_{\textnormal{\text{act}}}}(2\epsilon_1)^aD^k_0-\frac{\lambda_{v^0}\overline{D}_k}{2}>0}
		\end{align*}
		
		By construction, since we use new $(A^{k+1}_v(p))_{p\in \entier, v\in V}$ and $(D^{k+1}_m)_{m\in \entier}$ after the $k+1^{th}$ activation of $v^0$,  $((2\epsilon_1)^aD^k_{0}-\frac{\lambda_{v^0}\overline{D}_k}{2})_{k\in \entier}$ is an \textit{i.i.d.} family of random variables whose mean is given by
		\[\cesp{}{(2\epsilon_1)^aD^k_{v^0}-\frac{\lambda_{v^0}\overline{D}_k}{2}}=(2\epsilon_1)^a-\dfrac{\lambda_{v^0}\epsilon_0^a(2n-\lvert V'\rvert)}{2^{1+a}\lvert V'\rvert}.\]
		Notice that 
		\[\frac{\epsilon_0}{2}\left(\dfrac{\lambda_{v^0}(2n-\lvert V'\rvert)}{2^{1+a}\lvert V'\rvert}\right)^{1/a}>\frac{\epsilon_0}{2^{2+1/a}}\left(\min_v \lambda_v \right)^{1/a}>\epsilon_1,\]
		because $\lvert V'\rvert\leqslant n$.
		
		By the law of large numbers
		\[\cpro{}{\dfrac{1}{N^{1-a}K_{\textnormal{\text{act}}}}\sum_{k=1}^{N^{1-a}K_{\textnormal{\text{act}}}}\left((2\epsilon_1)^aD^k_{v^0}-\frac{\lambda_{v^0}\overline{D}_k}{2}\right)>0}\to 0 \text{ as } N\to +\infty,\]
		which proves the first part of the result.
		
		For any $t>0$, since by Lemma \ref{lem:notlong}
		\[\cpro{}{\inf_{v\in V', Ns\leqslant d(N^{1-a}K_{\textnormal{\text{act}}})}Q^N_v(s)>\frac{\epsilon_0}{2}}\to 1\text{ as }N\to +\infty,\] 
		we get as a consequence of the first part that 
		\begin{align}\label{eq:important}
		\cpro{}{\int_{0}^{d(N^{1-a}K)}\sigma_{v^0}(s)\d s> d(N^{1-a}K)\dfrac{\lambda_{v^0}}{2}, \sup_{Ns\leqslant d(N^{1-a}K)}Q^N_{v^0}(s)< \epsilon_1}\to 0 \text{ as } N\to +\infty.
		\end{align}
		With probability close to one, either 
		\[T^{\epsilon_1}_-(Q^N_{v^0})\leqslant d(N^{1-a}K_{\textnormal{\text{act}}}),\] 
		or 
		\[\int_{0}^{d(N^{1-a}K_{\textnormal{\text{act}}})}\sigma_{v^0}(s)\d s\leqslant  d(N^{1-a}K_{\textnormal{\text{act}}})\dfrac{\lambda_{v^0}}{2}.\]
		In the latter case, by \eqref{eq:traj},
		\[Q^N_{v^0}(\frac{b(N^{1-a}K_{\textnormal{\text{act}}})}{N})\geqslant \frac{d(N^{1-a}K_{\textnormal{\text{act}}})}{N}(\lambda_{v^0}-\frac{\lambda_{v^0}}{2}).\]
		More formally,
		\begin{multline*}
		\cpro{}{Q^N_{v^0}(\frac{d(N^{1-a}K_{\textnormal{\text{act}}})}{N})> d(N^{1-a}K_{\textnormal{\text{act}}})\dfrac{\lambda_{v^0}}{2}, \sup_{Ns\leqslant d(N^{1-a}K_{\textnormal{\text{act}}})}Q^N_{v^0}(s)< \epsilon_1}\\
		\leqslant	\cpro{}{\int_{0}^{d(N^{1-a}K_{\textnormal{\text{act}}})}\sigma_{v^0}(s)\d s> d(N^{1-a}K_{\textnormal{\text{act}}})\dfrac{\lambda_{v^0}}{2}, \sup_{Ns\leqslant d(N^{1-a}K_{\textnormal{\text{act}}})}Q^N_{v^0}(s)< \epsilon_1}\to 0 \text{ as } N\to +\infty
		\end{multline*}
			\qed

\textbf{{\underline{Proof of Lemma \ref{lem:nonneg}:}}}

		This proof relies on the coupling argument given in \eqref{eq:lowerd} (see the discussion there for the justification). In the event 
		\[\left\lbrace \inf_{v\in V', Ns\leqslant d(N^{1-a}K_{\textnormal{\text{act}}})}Q^N_v(s)\geqslant \frac{\epsilon_0}{2}\right\rbrace,\] it states that
		\begin{equation*}
		d(N^{1-a}K_{\textnormal{\text{act}}})\geqslant \sum_{k=1}^{N^{1-a}K_{\textnormal{\text{act}}}}\bar{D}_k= \sum_{k=1}^{N^{1-a}K_{\textnormal{\text{act}}}}\sum_{m=1}^{\check{G}(k)}(\frac{N\epsilon_0}{2})^a D^k_m,
		\end{equation*}
		with $\check{G}(k)$ \textit{i.i.d.} and independent from $(D^k_m)_{k,m\in \entier}$, of common distribution $\mathcal{G}(\frac{\lvert V'\rvert}{2n-\lvert V'\rvert})$. Obviously, because of that, 
		\begin{align*}
		\cpro{}{d(N^{1-a}K_{\textnormal{\text{act}}})<\frac{2N\epsilon_1}{\lambda_{v^0}},\; \inf_{v\in V', Ns\leqslant d(N^{1-a}K_{\textnormal{\text{act}}})}Q^N_v(s)\geqslant \frac{\epsilon_0}{2}}&\leqslant\cpro{}{ \sum_{k=1}^{N^{1-a}K_{\textnormal{\text{act}}}}\overline{D}_k<\frac{2N\epsilon_1}{\lambda_{v^0}}} \\
		&=\cpro{}{\frac{1}{N^{1-a}K_{\textnormal{\text{act}}}}\sum_{k=1}^{N^{1-a}K_{\textnormal{\text{act}}}}\overline{D}_k<\frac{2\epsilon_1}{K_{\textnormal{\text{act}}}\lambda_{v^0}}}
		\end{align*}
		By independence of the $\overline{D}_k$ and the law of large numbers, $\frac{1}{N^{1-a}K_{\textnormal{\text{act}}}}\sum_{k=1}^{N^{1-a}K_{\textnormal{\text{act}}}}\overline{D}_k\to \cesp{}{\overline{D}_1}$ almost surely. Elementary computations give
		\[\cesp{}{\overline{D}_k}=\dfrac{\epsilon_0^a(2n-\lvert V'\rvert)}{2^{a}\lvert V'\rvert},\]
		and thus since 
		\[K_{\textnormal{\text{act}}}= \dfrac{\epsilon_1K_{\textnormal{\text{time}}}}{n4^{1+a}\norm{q^0}{\infty}^a}=\epsilon_1\dfrac{1}{n^2\epsilon_0^a\norm{\lambda}{\infty}}.\]
		Notice that 
		\[\frac{2\epsilon_1}{\lambda_{v^0}K_{time}}=\frac{\min_v \lambda_v}{\lambda_{v^0}} \frac{\epsilon_0^a (2n-n)}{2^an}<\frac{\epsilon_0^a(2n-\lvert V' \rvert)}{2^a \lvert V'\rvert}= E[\overline{D}_k],\]
		
		where the last inequality is due to the fact that the right hand side is decreasing in $\lvert V'\rvert$.
		
		We get 
		\[\cpro{}{\dfrac{1}{N^{1-a}K_{\textnormal{\text{act}}}}\sum_{k=1}^{N^{1-a}K_{\textnormal{\text{act}}}}\overline{D}_k<\frac{2\epsilon_1}{K_{\textnormal{\text{act}}}\lambda_{v^0}}}\to 0 \text{ as } N\to +\infty.\]
		This entails
		\[	\cpro{}{d(N^{1-a}K_{\textnormal{\text{act}}})<\frac{2N\epsilon_1}{\lambda_{v^0}},\; E^N_-(\frac{d(N^{1-a}K_{\textnormal{\text{act}}}))}{N}}\to 0 \text{ as } N\to +\infty.\]
		By Lemma \ref{lem:notlong}, \[\cpro{}{E^N_-(\frac{d(N^{1-a}K_{\textnormal{\text{act}}}))}{N}}\to 1 \text{ as } N\to +\infty\]
		and the result is proved.	\qed 
	
	\section{\large{Conclusion and future research}}
	The main information to get from this article is that the homogenization result from \cite{thesis}\cite{us} is a good tool to study QB-CSMA as long as there is a separate method to ensure that all queues are positive. In particular, this suggest to look at a system of exponent $a_v$ such that for any $v\in V$, 
	\begin{equation}\label{eq:exta}
	\max_{\sigma\in S^*} \sigma_v\sum_{w\in V}a_w\sigma_w=c.
	\end{equation}
	This would ensure that every queue receives a non-zero service rate when it is of order $N$ in the fluid scale. Even when such a scheme is possible, in low load situation we can face an ODE which makes the reflection/absorption problem unavoidable.
	\subsection{\large{\underline{K-partite complete graphs}}}
	A natural class of graph to study in the vein of \eqref{eq:exta} are those where nodes can be disjointed between K part that interfere with each other. In other words, $V=\bigcup_{k=1,\ldots,p} V_p$, $p\leqslant n$ and $v_1\in V_k$, $v_2\in V_r$ such that $r\neq k$ implies $v_1\sim v_2$. For this type of interference graph, maximal stable sets do not intersect, each node will be in exactly one maximal stable set. We will once again have to distinguish cases: either
	
	\begin{itemize}
		\item there is a maximal stable set of size $\Upsilon'<\Upsilon$, or 
		\item all maximal stable sets have the same size $\Upsilon$.
	\end{itemize}
	In the first case, study of the ODE gives directly that some coordinates will reach zero before others because some nodes have a null service rate. The other case yields the same conclusion but requires more work: we present here the simplest case of a ``square interference graph''. This shows that there is little hope for functional limit results for general interference graph in low load situations using this method.
	
	Assuming homogenization, we have some ideas about the behaviour of the limit: the asymptotic service rates are given by 
	\[\bar{\pi}^q(1)=\bar{\pi}^q(3)=\dfrac{(q_1q_3)^a}{(q_1q_3)^a+(q_2q_4)^a},\, \bar{\pi}^q(2)=\bar{\pi}^q(4)=\dfrac{(q_2q_4)^a}{(q_1q_3)^a+(q_2q_4)^a}.\]
	Depending on arrival rates and initial condition, the solution to the ODE \eqref{eq:limit} can hit zero with one coordinate separately from others. In order to see this, notice that 
	\[(q^*_1-q^*_3)'(t)=\lambda_1-\lambda_3\text{ and } (q^*_2-q^*_4)'(t)=\lambda_2-\lambda_4.\]
	Notice as well that 
	\[s(q^*(t))=s(q^0)+(s(\lambda)-2)t\]
	up to the time a coordinate reaches 0. In fact the stability condition is $\lambda_1\vee\lambda_3+\lambda_2\vee\lambda_4<1$ so the expression for the sum will reach zero in finite time when the arrival rates are sub-critical.
	If 
	\[q^0_1-q^0_3=\lambda_1-\lambda_3=q^0_2-q_4^0=\lambda_2-\lambda_4=0,\]
	all coordinate will reach zero at the same time. Analysing the ODE for different parameter can yield different situations: \begin{itemize}
		\item if $\lambda_1<\lambda_3$ and $q^0_1\leqslant q^0_3$, the quantity $(q^*_3-q^*_1)(t)$is increasing from a positive value so it can never reach zero.
		\item if $\lambda_1<\lambda_3$ and $q^0_1>q^0_3$. There may exist a time such that $q^*_1(t)=q^*_3(t)=0$. Assuming the other coordinate do not reach zero before that point, it is only possible at time 
		\[t^*_1=\dfrac{q^0_1-q^0_3}{\lambda_3-\lambda_1}.\]
		We then have to do the same analysis with queues 2 and 4 to determine if it is possible to have all queues reach zero at the same time or not. In short, if $\lambda_2<\lambda_4$ and $q^0_2>q^0_4$ the condition to have all queues reach zero at the same time is 
		\[\dfrac{q^0_1-q^0_3}{\lambda_3-\lambda_1}=\dfrac{q^0_2-q^0_4}{\lambda_4-\lambda_2}=\dfrac{s(q^0)}{2-s(\lambda)}>0,\]
		otherwise, we have to compare the time it takes for $q^*_2$ and $q^*_4$ to become equal to the time it takes for $q^*_1$, $q^*_3$ and the sum to reach zero.
	\end{itemize}
	
	Although it is not proved that the total number of jobs in the system is monotonous in the arrival rates, a ``worst case scenario'' can be thought of as when $\lambda_1=\lambda_3$ and $\lambda_2=\lambda_4$. In this case, since 
	\[(q^*_1-q^*_3)'(t)=(q^*_2-q^*_4)'(t)=0,\]
	the two queues in a stable set cannot reach zero at the same time starting from different initial conditions. If the initial conditions are all equal, and the arrival rates are equal and subcritical, the solution of the ODE will reach zero in all coordinate simultaneously. In heavy load situation, the fluid limit does not reach 0 in any coordinate and it may become possible to analyse a faster time scale as in \cite{us}.
	\subsection{\large{\underline{Heavy Traffic}}}
	In the previous example, similarly to \cite{us}, for $a$ small enough, it may be possible to go further in time with time scales $N^{1+a}$, and $N^{1+2a}$ using homogenization if it is possible to bound queue lengths away from zero. Heuristically, in a super critical critical regime with $\lambda_1=\lambda_3$ and $\lambda_2=\lambda_4$, on the $N^{1+a}$ time scale, independent sets of size 1 will come into play. When all coordinates are bounded away from zero, any limiting point $\bar{q}$ should satisfy
	\[(\bar{q}_1-\bar{q}_3)'(t)=\dfrac{\bar{q}_3^a(t)-\bar{q}_1^a(t)}{(\bar{q}_1\bar{q}_3)^a(t)+(\bar{q}_2\bar{q}_4)^a(t)},\]
	and
	\[(\bar{q}_2-\bar{q}_4)'(t)=\dfrac{\bar{q}_4^a(t)-\bar{q}_2^a(t)}{(\bar{q}_1\bar{q}_3)^a(t)+(\bar{q}_2\bar{q}_4)^a(t)}.\]
	This would push queue lengths corresponding to queue in the same independent set towards equality.
	
	Similarly to \cite{us}, the time scale $N^{1+2a}$ will reveal the effect of idle time and we can expect some kind of state space collapse with $\widetilde{q}_1(t)=\widetilde{q}_3(t)$ and $\widetilde{q}_2(t)=\widetilde{q}_4(t)$ for any limiting point. The method by \cite{brams98} and \cite{will98} can most likely be adapted to that effect. A detailed analysis could potentially tell us that in heavy traffic for the square interference graph the delay scales like $\rho^{-1/2a}$ with $\rho $ the distance between arrival rates and the exterior of the stability region. 

	\appendix  \large
	\section{\large{Proof of Lemmas \ref{lem:uniquecig} and \ref{lem:qpos}, Existence and uniqueness of solutions to \eqref{eq:limitCIG}\label{app:ode}}}
\textbf{\underline{Proof of Lemma \ref{lem:uniquecig}:}}

		Let $q^N\in (0,+\infty)^V$ be such that $q^N \to q^0\in \reels_+^V$. Let the maximal solution to \eqref{eq:limit} on a complete interference graph with $\bar{q}^N(0)=q^N$ be denoted \[\bar{q}^N:D^N\to (0,+\infty)^V.\]
		
		For any $N$, there is $ T_{\text{ext}}(q^N)=\sup D^N>0$. Because of Lemma \ref{lem:qpos}, \[T_{\text{ext}}(q^N)=\tau^0(s\circ \bar{q}^N).\] Lemma \ref{lem:qpos} states that if $s(\lambda)>1$, $T_{\text{ext}}(q^N)=+\infty$, if $s(\lambda)<1$ $T_{\text{ext}}(q^N)=\frac{s(q^N)}{1-s(\lambda)}$ and if $s(\lambda)=1$, $T_{\text{ext}}(q^N)=+\infty$ but we will see that the limiting process may reach $\textbf{0}$ if $q^N\to \textbf{0}$.

		If $s(\lambda)< 1$ and $q^0=\textbf{0}$, $\tau^0(s\circ \bar{q}^N)\to 0$ as $N\to +\infty$ so $q^*(\cdot,q^N)\to \textbf{0}$ uniformly over compact time sets.
		
		If $s(\lambda)\geqslant 1$ or $q^0\neq \textbf{0}$, let us consider $D=(0,\liminf_N T_{\text{ext}}(q^N))$ and $\widetilde{q}^N: D\to \reels_+^V$, the restriction of $\bar{q}^N$ on $D$. The first step is to prove convergence of $\widetilde{q}^N$ for the uniform convergence and then extend to the positive half real line. Because $T_{\text{ext}}(q^N)=\frac{s(q^N)}{1-s(\lambda)}$ and $q^0\neq \textbf{0}$, we have $\liminf_N T_{\text{ext}}(q^N)>0$, and thus there is $t^0>0$ such that $(0,t^0)\subset D$. Since $q^N\to q^0\in \reels_+^V$, $\norm{g}{\infty}\leqslant 1$ and $(\widetilde{q}^N)'(t)=g(\widetilde{q}^N)(t)$,  $\widetilde{q}^N$ is relatively compact in $C(I,\reels_+^V)$ for the convergence uniform over compact time sets by the Arzelà-Ascoli theorem (see Theorem 7.2 in Chapter 2 of \cite{bill}). Indeed,  for any $\delta>0$,
		\[\sup_{ t,s\leqslant T, \lvert t-s\rvert\leqslant \delta}\norm{q^*_v(t)-q^*_v(s)}{\infty}\leqslant \norm{g}{\infty}\delta .\]
		Let $N_k$ be a subsequence such that  $\widetilde{q}^{N_k}\to\bar{q}$ uniformly over compact sets of $D$. Note that compact sets of $D$ are bounded away from $0$. We will prove that $\bar{q}$ must be a solution to \eqref{eq:limitCIG}.

		First we prove that $g(\widetilde{q}^N)\to g(\bar{q})$ uniformly over compact time sets of $D$. Fix a time horizon $T<\sup I$ and $\nu>0$. By definition of $D$, 
		\[\inf_{\nu\leqslant t\leqslant T, v\in V}\widetilde{q}^N_v(t)\geqslant \eta>0.\] First, $g$ is Lipschitz continuous on $(\eta,C_T)^V$ because it has bounded derivative, so it is uniformly continuous. By definition of uniform continuity, for any $\epsilon>0$, there exists $\delta>0$ such that for any $q,q'\in (\eta,C_T)^V$ with $\norm{q-q'}{\infty}\leqslant \delta$,
		\[\norm{g(q)-g(q')}{\infty}\leqslant \epsilon.\]
		Let $N_0$ be such that $\sup_{t\leqslant T}\norm{\widetilde{q}^N(t)-\bar{q}(t)}{\infty}\leqslant \delta$ for any $N\geqslant N_0$. Then for any $N>N_0$, 
		\[\sup_{0<t \leqslant T}\norm{g(\widetilde{q}^N)(t)-g(\bar{q}(t))}{\infty}\leqslant \epsilon,\]
		i.e. $g(\widetilde{q}^N)\to g(\bar{q})$ uniformly over compact time sets of $D$.
		
		Since $q^N\to q^0\in \reels^V_+$, $(\widetilde{q}^N)'=g(\widetilde{q}^N)$, $\widetilde{q}^N\to \bar{q}$ and $g(\widetilde{q}^N)\to g(\bar{q})$ uniformly over compact time sets of $(0,t^0)$,  by Theorem 7.17 of \cite{rud}, $\bar{q}$ is differentiable and $(\widetilde{q}^N)'\to \bar{q}'$ uniformly over compact time sets of $(0,t^0)$. For any $0<t< t^0$, by uniqueness of the limit $\bar{q}'(t)=g(\bar{q})(t)$, and thus $\bar{q}$ must be a continuous function such that for any $t\in(0,t^0)$ we have $\bar{q}'(t)=g(\bar{q})(t)$. This also proves the existence of solutions to \eqref{eq:limitCIG}. We only used Lemma \ref{lem:qpos} for solutions of \eqref{eq:limit} with positive initial condition so there is no logic loop.
		
		We now prove uniqueness of solutions to \eqref{eq:limitCIG}. Let  $\overline{q}:D_1\to \reels_+^V$ and $\widetilde{q}:D_2\to \reels_+^V$ be two solutions starting from $ q^0$ at time $t=0$. For $t\in D_1\cap D_2$, the sum can be expressed as
		
		\[\zeta(t)\coloneqq \max(s(q^0)+(s(\lambda)-1)t,0).\]
		
		If $s(\lambda)\geqslant 1$ and $q^0\neq \textbf{0}$, $\zeta$ will never reach $0$. If $s(\lambda)< 1 $, it will in finite time.

		Let $t_d\in D_1\cap D_2$ be the first time that $\overline{q}(t)\neq \widetilde{q}(t)$.	We first make the assumption that $t_d$ is also the infimum of times such that $s(\overline{q}(t)^a)> s(\widetilde{q}^a(t))$ and deal with the complementary later. By continuity there is $t_d>0$ and $\delta>0$ such that $\overline{q}(s)\neq \widetilde{q}(s)$ and $s(\overline{q}^a(s))> s(\widetilde{q}^a(s))$for any  $s\in (t_d,t_d+\delta)\subset D_1\cap D_2$. 
		
		For any $v$ such that $\overline{q}_v(s)< \widetilde{q}_v(s)$, we have $\frac{\overline{q}_v(s)}{s(\overline{q}^a(s))}< \frac{\widetilde{q}_v(s)}{s(\widetilde{q}^a(s))}$ and thus, 
		
		\[\overline{q}_v(s)- \widetilde{q}_v(s)< 0 \Rightarrow (\overline{q}_v'(s)- \widetilde{q}_v'(s))> 0.\]
		
		By continuity of the solutions, since $\bar{q}(t_d)=\widetilde{q}(t_d)$, $s(\overline{q}^a(s))> s(\widetilde{q}^a(s))$ implies $\overline{q}_v(s)\geqslant \widetilde{q}_v(s)$ for all $v\in V$, for all $s\in (t_d,t_d+\delta)$ because as soon as $\bar{q}_v(t)$ would become smaller than $\widetilde{q}_v(t)$, $\bar{q}_v(t)$ would become increasing. Since $s(\overline{q})=s(\widetilde{q})$ the only possibility is that  $\overline{q}_v(s)= \widetilde{q}_v(s)$ for all $v\in V$ which contradicts $s(\overline{q}(s)^a)> s(\widetilde{q}^a(s))$ for any  $s\in (t_d,t_d+\delta)$. The first time $\widetilde{q}$ and $\bar{q}$ are different is different from the first time $s\circ \widetilde{q}^a$ and $s\circ\bar{q}^a$ are different.
		
		Let us introduce for any $t\in [t_d,t_d+\delta)$, 
		
		\[\zeta_a(t)\coloneqq s(\overline{q}^a(t))=s(\widetilde{q}^a(t)).\]
		
		We can rewrite \eqref{eq:limitCIG} as
		
		\begin{equation}
		\left\lbrace\begin{array}{ll}
		f'(t)&=\lambda-\dfrac{f(t)^a}{\zeta_a(t)} \text{ if $\zeta_a(t)>0$} \\
		s(f)&=\zeta\\
		f(0)&=q^0 \end{array}\right. \label{eq:limitCIGfinal}
		\end{equation}
		
		For any $v\in V$, and $t>0$, such that $\overline{q}_v(t)- \widetilde{q}_v(t)\geqslant 0$ we have 
		
		\[\overline{q}_v'(t)- \widetilde{q}_v'(t)=\frac{\widetilde{q}_v(t)^a-\overline{q}_v(t)^a}{\zeta_a(t)}\leqslant 0\text{ because }a>0.\]
		
		Once again by continuity of the solutions, this implies $\overline{q}_v(t)\leqslant \widetilde{q}_v(t)$ for any $t\in D_1\cap D_2$ and equality since $s(\overline{q})=s(\widetilde{q})$. Thus there is a unique maximal solution defined up to the exit time of $(0,+\infty)^V$. By Lemma \ref{lem:qpos}, \[\lim_{t\to T_{\text{ext}(q^0)}}\bar{q}(t)=\lim_{t\to T_{\text{ext}(q^0)}}\widetilde{q}(t)=\textbf{0}.\] We extend the solution as in Definition \ref{def:ext} to obtain convergence on $\reels_+$. 
	\qed 
	
	We now turn to the proof of Lemma \ref{lem:qpos}.\vspace{0.3cm}
	
\textbf{{\underline{Proof of Lemma \ref{lem:qpos}:}}}
		Let $q^0\in \reels_+^V$. If $q^0\in (0,+\infty)^V$, let $q^*=q^*(\cdot,q^0)$ solution to \eqref{eq:limit}. If $q^0\in \reels_+^V\setminus(0,+\infty)^V$, assume the existence of a solution (they necessarily exist for $q^0\in \reels_+^V\setminus\lbrace\textbf{0}\rbrace$, see the proof of Lemma \ref{lem:uniquecig} for more details). We use $q^*$ to denote an arbitrary solution to \eqref{eq:limitCIG}. If $q^0\neq\textbf{0}$, the expression of $\tau^0(s\circ q^*)$ is a direct consequence of the expression of the sum.
		
		For any $q\in \entier^V$, there always exists $q_v\geqslant \frac{s(q)}{n}$ because otherwise $\frac{ns(q)}{n}<s(q)$ by summing over $V$. Hence, $s(q^a)\geqslant (\frac{s(q)}{n})^a$. 
		
		We now prove that $\tau^0(s\circ q^*)=\tau^0(q^*)$ if $s(\lambda)>1$. For any $q^0\in \reels_+^V$, $t>0$ we get
		\begin{align*}(q^*_v)'(t)&=\lambda_v-\dfrac{(q^*_v(t))^a}{s((q^*(t))^a)}\\
		&\geqslant \lambda_v-\left(\dfrac{nq^*_v(t)}{s(q^*(t))}\right)^a\\
		&= \lambda_v-\left(\dfrac{nq^*_v(t)}{(s(\lambda)-1)t}\right)^a.
		\end{align*}
		
		If $q^*_v(t)< (\frac{\lambda_v}{2})^{1/a} \frac{(s(\lambda)-1)}{n}t$ we get $(q^*_v)'(t)> \frac{\lambda_v}{2}$.	If $\frac{\lambda_v}{2}< (\frac{\lambda_v}{2})^{1/a} \frac{(s(\lambda)-1)}{n}$, $q^*_v(t)$ cannot become smaller than $\frac{\lambda_v}{2}t$: if at some point $q^*_v(t)\leqslant(\frac{\lambda_v}{2})^{1/a} \frac{(s(\lambda)-1)}{n}t$, by continuity of $q^*$, for any $s>0$, 
		\[q^*_v(t+s)\geqslant q^*_v(t)+\frac{\lambda_v}{2}s.\]
		
		Similarly,  If $\frac{\lambda_v}{2}\geqslant (\frac{\lambda_v}{2})^{1/a} \frac{(s(\lambda)-1)}{n}$, by continuity, $q^*_v(t)$ cannot become smaller than $(\frac{\lambda_v}{2})^{1/a} \frac{(s(\lambda)-1)}{n}t$  because as soon as $q^*_v(t)\leqslant (\frac{\lambda_v}{2})^{1/a} \frac{(s(\lambda)-1)}{n}t$, we get \[(q^*_v)'(t)\geqslant \frac{\lambda_v}{2}\geqslant(\frac{\lambda_v}{2})^{1/a} \frac{(s(\lambda)-1)}{n} .\] 
		
		Thus for every $v\in V$, and $t\geqslant 0$, 
		
		\[q^*_v(t)\geqslant \min\left(\frac{\lambda_v}{2}, (\frac{\lambda_v}{2})^{1/a} \frac{(s(\lambda)-1)}{n}\right)t,\] 
		which implies $\tau^0(q^*)=\tau^0(s\circ q^*)=+\infty$. 
		
		We next prove that $\tau^0(q^*)= \tau^0(s\circ q^*)$ if $q^0\neq\textbf{0}$ and $s(\lambda)\leqslant 1$. The first step is done now: we prove that coordinates which started empty become positive for any small enough positive time. If $s(\lambda)=1$, $s(q^*)$ is constant and thus always greater than $\frac{s(q^0)}{2}$. If $s(\lambda)<1$, for $t\leqslant \frac{s(q^0)}{2\lvert s(\lambda)-1\rvert}$, we have $s(q^*(t))\geqslant \frac{s(q^0)}{2}$. In both cases, there is $t^0>0$ and $\epsilon>0$ such that $\inf s(q^*)(t)\geqslant \epsilon$ for any $t\leqslant t^0$. Moreover, for any $\epsilon>0$, if $s(q^*(t))\geqslant \epsilon$, we can give a lower bound to the derivative in time for solutions to \eqref{eq:limitCIG}. Using the definition of the ODE, and $s(q^a)\geqslant \left(\frac{s(q)}{n}\right)^a$, we get 
		\begin{align}(q^*_v)'(t)&=\lambda_v-\dfrac{(q^*_v(t))^a}{s((q^*(t))^a)}\notag\\
		&\geqslant \lambda_v-\left(\dfrac{nq^*_v(t)}{s(q^*(t))}\right)^a\label{eq:derivv}\\
		&\geqslant \lambda_v-\left(\dfrac{nq^*_v(t)}{\epsilon}\right)^a\notag
		\end{align}
		For any $t\geqslant 0$ if $s(\lambda)=1$ or $t\leqslant\frac{s(q^0)}{2\lvert s(\lambda)-1\rvert}$ if $s(\lambda)<1$, whenever $q^*_v(t)\leqslant\left(\frac{\lambda_v}{2}\right)^{1/a}\frac{s(q^0)}{2n} $, it follows that $(q^*_v)'(t)>\frac{\lambda_v}{2}$. So there is $t^0,\epsilon'>0$ such that any $v$ with $q^0_v\leqslant \epsilon'$ gets $q^*_v$ increasing for $t\in (0,t^0)$. Even if some queues start null, they become positive for $t>0$ small enough.
		
		We now prove that if they do touch $0$, all coordinates reach it simultaneously. It is only possible for queue lengths to touch 0 if $s(\lambda)<1$. By \eqref{eq:derivv}, if $q^*_v(t)\leqslant \left(\frac{\lambda_v}{2}\right)^{1/a}\frac{s(q^*(t))}{n}$, we have $(q^*_v)'(t)\geqslant \frac{\lambda_v}{2}$. If $\tau^0(q^*)<\tau^0(s\circ q^*)$, this would mean that $s(q^*)(\tau^0(q^*))>0$. Without loss in generality, assume that $q^*_v(t)\to 0$ as $ t\to \tau^0(q^*)$. Since $q^*_v$ is continuous and $q^*_v(t)\to 0$ as $t\to\tau^0(q^*)$, there exists a time interval $(t_-, \tau^0(q^*))$ such that $q^*_v(t)\leqslant \left(\frac{\lambda_v}{2}\right)^{1/a}\frac{s(q^*)(\tau^0(q^*))}{n}$ and $q^*_v$ is decreasing. Since $s(\lambda)<1$, $s\circ q^*$ is decreasing so $s( q^*)(t)\geqslant s(q^*)(\tau^0(q^*))$ for any $t\in (t_-, \tau^0(q^*))$. This leads to a contradiction because $(q^*_v)'(t)\geqslant \frac{\lambda_v}{2}$ when $q^*_v(t)\leqslant \left(\frac{\lambda_v}{2}\right)^{1/a}\frac{s(q^*(t))}{n}$. Hence $\tau^0(q^*)=\tau^0(s\circ q^*)$. 	\qed 
	
\thanks{The author would like to thank Florian Simatos and Laurent Miclo for the numerous discussions on the subject that provided moral support and concrete leads to overcome some of the difficulties treated in this paper.}

	%
	%
	%
	%
	%
	%
	%
	
	\bibliographystyle{plainurl}
	
\bibliography{preprintfluid}


\end{document}